\documentclass[stsy]{informs-stsy}

\usepackage{epsf}
\usepackage{graphicx}
\usepackage{amsfonts,amssymb,amsmath}
\usepackage{algorithmic,algorithm}
\usepackage{subfigure}
\usepackage{latexsym}
\usepackage{bm}
\usepackage{url}

\usepackage{natbib}
 \bibpunct[, ]{(}{)}{,}{a}{}{,}%
 \def\bibfont{\small}%
\usepackage[normalem]{ulem}

\newcommand{\expect}[1]{\mathbb{E}\left[#1\right]}

\newtheorem{theorem}{\textbf{Theorem}}
\newtheorem{lemma}{\textbf{Lemma}}
\newtheorem{corollary}{\textbf{Corollary}}
\newtheorem{assumption}{\textbf{Assumption}}
\newtheorem{definition}{\textbf{Definition}}

\newtheorem{remark}{\textbf{Remark}}
\newtheorem{prop}{Proposition}

\usepackage[colorlinks=true,breaklinks=true,bookmarks=true,urlcolor=blue,
     citecolor=blue,linkcolor=blue,bookmarksopen=false,draft=false]{hyperref}

\renewcommand{\theARTICLETOP}{}

\begin{document}
 \RUNAUTHOR{X. Wei and M. J. Neely}
\RUNTITLE{Opportunistic Scheduling over Renewal Systems}
\TITLE{Opportunistic Scheduling over Renewal Systems: An Empirical Method}

\ARTICLEAUTHORS{%
\AUTHOR{Xiaohan Wei}
\AFF{Department of Electrical Engineering, University of Southern California \\
\EMAIL{xiaohanw@usc.edu}}
\AUTHOR{Michael J. Neely}
\AFF{Department of Electrical Engineering, University of Southern California \\
\EMAIL{mjneely@usc.edu}}
} 

\ABSTRACT{
This paper considers an opportunistic scheduling problem over a renewal system. A controller observes a 
random event at the beginning of each renewal frame and then chooses an action in response to the event, which affects the duration of the frame, the amount of resources used, and a penalty metric. The goal is to make frame-wise decisions so as to minimize the time average penalty subject to time average resource constraints. This problem has applications to task processing and communication in data networks, as well as to certain classes 
of Markov decision problems. 
We formulate the problem as a dynamic fractional program and propose an adaptive algorithm which uses an empirical accumulation as a feedback parameter. A key feature of the proposed algorithm is that it does not require knowledge of the random event statistics and potentially allows (uncountably) infinite event sets. 
We prove the algorithm satisfies all desired constraints and achieves $O(\epsilon)$ near optimality with probability 1. 
}


\MSCCLASS{90C15, 90C30, 90C40, 93E35}
\HISTORY{}

\maketitle

\section{Introduction}
\begin{figure}[htbp]
   \centering
   \includegraphics[height=1in]{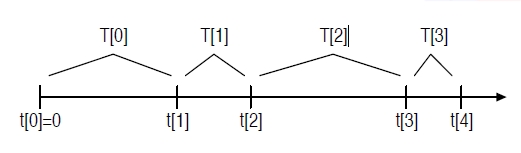} 
   \caption{An illustration of a sequence of renewal frames.}
   \label{fig:renewal}
\end{figure}

Consider a system that operates over the timeline of real numbers $t \geq 0$.  The timeline is divided into back-to-back periods called \emph{renewal frames} and the start of each frame is called a \emph{renewal} (see Fig. \ref{fig:renewal}).   The system state is refreshed at each renewal.    At the start of each renewal frame $n \in \{0, 1, 2, \dots\}$ the controller observes a random event $\omega[n]\in\Omega$ and then takes an action $\alpha[n]$ from an action set $\mathcal{A}$ in response to $\omega[n]$.   The pair $(\omega[n], \alpha[n])$ affects: (i)  the duration of that renewal frame; (ii)  a vector of resource expenditures for that frame;  (iii) a penalty incurred on that frame.   The goal is to choose actions over time to minimize time average penalty subject to time average constraints on the resources without knowing any statistic of $\omega[n]$. We call such a problem \textit{opportunistic scheduling over renewal systems}. 

\subsection{Example applications} 

This problem has applications to task processing in computer networks, 
and certain generalizations of Markov decision problems.  

\begin{itemize} 
\item Task processing networks:  Consider a device that processes tasks back-to-back.  Each renewal period corresponds to the time required to complete a single task.  The random event $\omega[n]$ observed corresponds to a vector of task parameters, including the type, size, and resource requirements for that particular task.  The action
set $\mathcal{A}$ consists of different processing mode options, and the specific action  $\alpha[n]$ determines the processing time, energy expenditure, and task quality.  In this case, task quality can be defined as a negative penalty, and the goal is to maximize time average quality subject to power constraints and task completion rate constraints. A specific example of this sort is the following file downloading problem: Consider a wireless device that repeatedly downloads files.  The device has two states: \emph{active} (wants to download a file) and \emph{idle} (does not want to download a file). Renewals occur at the start of each new active state.
Here, $\omega[n]$ denotes the observed wireless channel state, which affects the success probability of downloading a file (and thereby affects the transition probability from active to idle).  This example is discussed further in the simulation section (Section \ref{simulation}). 

\item Hierarchical Markov decision problems:  Consider a slotted \emph{two-timescale Markov decision processes (MDP)}
over an infinite horizon and with constraints on average cost per slot.  
An MDP is run on the lower level, with a special state that is recurrent under any sequence of actions. The
renewals are defined as revisitation times to that state.  On a higher level, 
a random event $\omega$ is observed upon each revisitation to the renewal state on the lower level. Then, a decision is made on the higher level in response to $\omega$, which in turn affects the transition probability and penalty/cost received per slot on the lower level until the next renewal.
Such a problem is a generalization of classical MDP problem (e.g. \citet{Ro02}, \citet{Be01}) and has been considered previously in \citet{wernz2013multi}, \citet{chang2003multitime} with discrete finite state and full information on both levels. A heuristic method is also proposed in  \citet{wernz2013multi} when some of the information is unknown. The algorithm of the current paper does not require 
knowledge of the statistics of $\omega$ and allows the event set $\Omega$ to be potentially (uncountably) infinite.
\end{itemize} 

\subsection{Previous approaches on renewal systems} 
Most works on optimization over renewal systems consider the simpler scenario of knowing the probability distribution of $\omega[n]$. In such a case, one can show via the renewal-reward theory that the problem can be solved (offline) by finding the solution to a linear fractional program. This idea has been applied to solve MDPs in the seminal work \citet{Fo66}. 
Methods for solving linear fractional programs can also be found, for example, in \citet{Sc83, BV04}. 
However, the practical limitations of such an offline algorithm are twofold: First, if the event set 
$\Omega$ is large, then, there are too many probabilities $Pr(\omega[n] = \omega),~\omega\in\Omega$ to estimate and the corresponding offline optimization problem may be difficult to solve even if all probabilities are estimated accurately. Second, generic offline 
optimization solvers may not take advantage of the special renewal 
structure of the system.  One notable example is the treatment of 
power and delay minimization for a multi-class M/G/1 queue in \citet{Yao02, LN14}, 
where the renewal structure allows a well known $c$-$\mu$ rule 
for delay minimization to be extended to treat both power and 
delay constraints.


The work in \citet{Neely2010,Ne09} presents a new \emph{drift-plus-penalty (DPP) ratio} algorithm solving renewal optimizations knowing the distribution of $\omega[n]$. The algorithm treats the constraints via \emph{virtual queues} so that one only requires to minimize an unconstrained ratio during every renewal frame. The algorithm provably meets all constraints and achieves asymptotic
near-optimality. The works  \citet{wang2015dynamic, urgaonkar2015dynamic} show that the edge cloud server migration problem can be formulated as a specific renewal optimization. Using a variant of the DPP ratio algorithm, they show that solving a simple stochastic shortest path problem during every renewal frame gives near-optimal performance. The work \citet{wei2018asynchronous} solves a more general asynchronous optimization over parallel renewal systems, though the knowledge of the random event statistics is still required. It is worth noting that the work \citet{Ne09} also proposes a heuristic algorithm when the distribution of 
$\omega[n]$ is not known. That algorithm is partially analyzed:  It is shown that if a certain process converges, then the algorithm converges to a near-optimal point.  However, whether or not such a process converges is unknown.


\subsection{Other related works}
The renewal optimization problem considered in this paper is a generalization of stochastic optimization
over fixed time slots. Such problems are categorized based on whether or not the random event is observed
before the decision is made.  Cases where the random event is observed before taking actions are often referred to as \textit{opportunistic scheduling problems}. Over the past decades, many algorithms have been proposed including max-weight (\citet{tassiulas1990stability, tassiulas1993dynamic}), Lyapunov optimization (\citet{eryilmaz2006joint, eryilmaz2007fair, Neely2010, georgiadis2006resource}), fluid model methods (\citet{stolyar2005maximizing, eryilmaz2007fair}), and dual subgradient methods (\citet{lin2004joint, ribeiro2010ergodic}) are often used.  

Cases where the random events are not observed are referred to as  \textit{online learning problems}. Various algorithms are developed for unconstrained learning including the weighted majority algorithm (\citet{littlestone1994weighted}), multiplicative weighting algorithm (\citet{freund1999adaptive}), following the perturbed leader (\citet{hutter2005adaptive}) and online gradient descent (\citet{zinkevich2003online, hazan2014beyond}). The resource constrained learning problem is studied in \citet{mahdavi2012trading} and \citet{wu2015algorithms}.
Online learning with an underlying MDP structure is also treated using modified multiplicative weighting (\citet{even2005experts}) and improved following the perturbed leader (\citet{yu2009markov}).



\subsection{Our contributions} 
In this work, we focus on opportunistic scheduling over renewal systems and propose a new algorithm that runs online (i.e. takes actions in response to each observed 
$\omega[n]$). Unlike prior works, the proposed algorithm requires neither the statistics of $\omega[n]$ nor explicit estimation of them, and is fully analyzed with convergence properties that 
hold with probability 1. From a technical perspective, we prove near-optimality of the algorithm by showing  
asymptotic stability of a customized process, relying on a novel construction of exponential supermartingales which could be of independent interest.
We complement our theoretical results with
simulation experiments on a time varying constrained MDP.

\section{Problem Formulation and Preliminaries}\label{formulation}
Consider a system where the time line is divided into back-to-back time periods called frames. At the beginning of frame $n$ ($n\in\{0,1,2,\cdots\}$), a controller observes the realization of a random variable $\omega[n]$, which is an i.i.d. copy of a random variable taking values in a compact set $\Omega\in\mathbb{R}^q$ with distribution function unknown to the controller.
Then, after observing the random event $\omega[n]$, the controller chooses an action vector $\alpha[n]\in\mathcal{A}$. Then, the tuple $(\omega[n],~\alpha[n])$ induces the following random variables:
\begin{itemize}
\item The penalty received during frame $n$: $y[n]$.
\item The length of frame $n$: $T[n]$.
\item A vector of resource consumptions during frame $n$:
$\mathbf{z}[n]=[z_1[n],~z_2[n],~\cdots,~z_L[n]]$.
\end{itemize}
We assume that \textit{given $\alpha[n]=\alpha$ and $\omega[n]=\omega$ at frame $n$, $(y[n],T[n],\mathbf{z}[n])$ is a random vector independent of the outcomes of previous frames}, with \emph{known} expectations. We then denote these conditional expectations as 
\begin{align*}
\hat{y}(\omega,\alpha)=&\expect{y[n]~|~\omega,\alpha},\\
\hat{T}(\omega,\alpha)=&\expect{T[n]~|~\omega,\alpha},\\
\hat{\mathbf{z}}(\omega,\alpha)=&\expect{\hat{\mathbf{z}}[n]~|~\omega,\alpha},
\end{align*}
which are all deterministic functions of $\omega$ and $\alpha$. This notation is useful when we want to highlight the action $\alpha$ we choose. 
The analysis assumes a single action in response to the observed $\omega[n]$ at each frame. Nevertheless, an ergodic MDP can fit into this model by defining the action as a selection of a policy to implement over that frame so that the corresponding $\hat{y}(\omega,\alpha)$, $\hat{T}(\omega,\alpha)$ and $\hat{\mathbf{z}}(\omega,\alpha)$ are expectations over the frame under the chosen policy.

Let
\begin{align*}
\overline{y}[N]&=\frac1N\sum_{n=0}^{N-1}y[n],\\
\overline{T}[N]&=\frac1N\sum_{n=0}^{N-1}T[n],\\
\overline{z}_l[N]&=\frac1N\sum_{n=0}^{N-1}z_l[n]~~~l\in\{1,2,\cdots,L\}.
\end{align*}
The goal is to minimize the time average penalty subject to $L$ constraints on resource consumptions. Specifically, we aim to solve the following fractional programming problem:
\begin{align}
\min~~&\limsup_{N\rightarrow\infty}\frac{\overline{y}[N]}{\overline{T}[N]}\label{prob-1}\\
\textrm{s.t.}~~& \limsup_{N\rightarrow\infty}\frac{\overline{z}_l[N]}{\overline{T}[N]}\leq c_l,~~\forall l\in\{1,2,\cdots,L\},\\
&\alpha[n]\in\mathcal{A},~\forall n\in\{0,1,2,\cdots\} \label{prob-3},
\end{align}
where $c_l,~l\in\{1,2,\cdots,L\}$ are nonnegative constants, and both the minimum and constraint are taken in an almost sure sense. Finally, we use $\theta^*$ to denote the minimum that can be achieved by solving above optimization problem. For simplicity of notation, let
\begin{equation}\label{def-K}
K[n]=\sqrt{\sum_{l=1}^L(z_l[n]-c_lT[n])^2}.
\end{equation}

\subsection{Assumptions}
Our main result requires the following assumptions, their importance will become clear as we proceed. We begin with the following boundedness assumption:
\begin{assumption}[Exponential type]\label{bounded-assumption}
Given $\omega[n]=\omega\in\Omega$ and $\alpha[n]=\alpha\in\mathcal{A}$ for a fixed $n$, it holds that 
$T[n]\geq1$ with probability 1 and $y[n],~K[n],~T[n]$ are of exponential type, i.e. there exists a constant $\eta>0$ s.t.
\begin{align*}
&\expect{\left.\exp\left(\eta \big|y[n]\big|\right)~\right|\omega,\alpha}\leq B+1,\\
&\expect{\left.\exp\left(\eta \big|K[n]\big|\right)~\right|\omega,\alpha}\leq B+1,\\
&\expect{\left.\exp\left(\eta \big|T[n]\big|\right)~\right|\omega,\alpha}\leq B+1,
\end{align*}
where $B$ is a positive constant.
\end{assumption}

The following proposition is a simple consequence of the above assumption:
\begin{prop}\label{prop-1}
~~Suppose Assumption \ref{bounded-assumption} holds.
Let $X[n]$ be any of the three random variables $y[n]$, $K[n]$ and $T[n]$ for a fixed $n$. Then, given 
any $\omega[n]=\omega\in\Omega$ and $\alpha[n]=\alpha\in\mathcal{A}$,
\begin{align*}
\expect{\left.\big|X[n]\big|~\right|\omega,\alpha}\leq B/\eta,~~
\expect{\left.X[n]^2~\right|\omega,\alpha}\leq 2B/\eta^2.
\end{align*}
\end{prop}
The proof follows from the inequality: 
$$
B+1\geq\expect{\left.e^{\eta \big|X[n]\big|}~\right|\omega,\alpha}
\geq 1+\eta\cdot \expect{\left.\big|X[n]\big|~\right|\omega,\alpha}
+ \frac{\eta^2}{2}\cdot\expect{\left.X[n]^2~\right|\omega,\alpha}.
$$ 

\begin{assumption}\label{optimal-assumption}
~~There exists a positive constant $\theta_{\max}$ large enough so that the optimal objective of $\eqref{prob-1}-\eqref{prob-3}$, denoted as $\theta^*$, falls into $[0,\theta_{\max})$ with probability 1.
\end{assumption}

\begin{remark}
~~If $\theta^*<0$, then, we shall find a constant $c$ large enough so that $\theta^*+c\geq0$. Then, define a new penalty $y'[n]=y[n]+cT[n]$. It is easy to see that minimizing $\limsup_{N\rightarrow\infty}\overline{y}[N]/\overline{T}[N]$ is equivalent to minimizing $\limsup_{N\rightarrow\infty}\overline{y'}[N]/\overline{T}[N]$ and the optimal objective of the new problem is $\theta^*+c$, which is nonnegative. 
\end{remark}

\begin{assumption}\label{assumption-for-algorithm}
~~Let $\left(\hat y(\omega,\alpha),~\hat T(\omega,\alpha),~\hat{\mathbf{z}}(\omega,\alpha)\right)$ be the performance vector under a certain $(\omega,\alpha)$ pair. Then, for any fixed $\omega\in\Omega$, the set of achievable performance vectors over all $\alpha\in\mathcal{A}$ is compact.
\end{assumption}

In order to state the next assumption, we need the notion of  \textit{randomized stationary policy}. We start with the definition:
\begin{definition}[Randomized stationary policy]\label{RSP}
A randomized stationary policy is an algorithm that at the beginning of each frame $n$, after observing the random event $\omega[n]$, the controller chooses $\alpha^*[n]$ with a conditional probability that is the same for all $n$. 
\end{definition}

\begin{assumption}[Bounded achievable region]\label{compact}
Let
\[(\overline{y},~\overline{T},~\overline{\mathbf{z}})\triangleq\expect{(\hat y(\omega[0],\alpha^*[0]),~\hat T(\omega[0],\alpha^*[0]),~\hat{\mathbf{z}}(\omega[0],\alpha^*[0]))}\]
be the one-shot average of one randomized stationary policy. Let $\mathcal{R}\subseteq\mathbb{R}^{L+2}$ be the set of all achievable one-shot averages $(\overline{y},~\overline{T},~\overline{\mathbf{z}})$. Then, $\mathcal{R}$ is bounded.
\end{assumption}

\begin{assumption}[$\xi$-slackness]\label{slack}
 There exists a randomized stationary policy $\alpha^{(\xi)}[n]$ such that the following holds,
\[\frac{\expect{\hat{z}_l\left(\omega[n],\alpha^{(\xi)}[n]\right)}}
{\expect{\hat{T}(\omega[n],\alpha^{(\xi)}[n])}}=c_l-\xi,~~\forall l\in\{1,2,\cdots,L\},\]
where $\xi>0$ is a constant.
\end{assumption}

\begin{remark}[Measurability issue]
~~We implicitly assume the policies for choosing $\alpha$ in reaction to $\omega$ result in a measurable $\alpha$, so that $T[n]$, $y[n]$, $\mathbf{z}[n]$ are valid random variables and the expectations in Assumption \ref{compact} and \ref{slack} are well defined. This assumption is mild. For example, when the sets $\Omega$ and $\mathcal{A}$ are finite, it holds for any randomized stationary policy. More generally, if $\Omega$ and $\mathcal{A}$ are measurable subsets of some separable metric spaces, this holds whenever the conditional probability in Definition \ref{RSP} is ``regular'' (see \citet{Durrett} for discussions on regular conditional probability), and $T[n]$, $y[n]$, $\mathbf{z}[n]$ are continuous functions on $\Omega\times\mathcal{A}$. 
\end{remark}

\section{An Online Algorithm}\label{online-section}
We define a vector of virtual queues $\mathbf{Q}[n]=[Q_1[n]~Q_2[n]~\cdots~Q_L[n]]$ which are 0 at $n=0$ and updated as follows:
\begin{equation}\label{queue-update}
Q_l[n+1]=\max\{Q_l[n]+z_l[n]-c_lT[n],0\}.
\end{equation}
The intuition behind this virtual queue idea is that if the algorithm can stabilize $Q_l[n]$, then the ``arrival rate'' $\overline{z}_l[N]/\overline{T}[N]$ is below ``service rate'' $c_l$ and the constraint is satisfied. The proposed algorithm then proceeds as in Algorithm \ref{online-algorithm} via two fixed parameters $V>0$, $\delta>0$, and an additional process $\theta[n]$ that is initialized to be $\theta[0]=0$.
For any real number $x$, the notation $[x]_0^{\theta_{\max}}$ stands for ceil and floor function:
\[[x]_0^{\theta_{\max}}=\left\{
                     \begin{array}{ll}
                       \theta_{\max}, & \hbox{if $x\in(\theta_{\max},+\infty)$;} \\
                       x, & \hbox{if $x\in[0,\theta_{\max}]$;} \\
                       0, & \hbox{if $x\in(-\infty, 0)$.}
                     \end{array}
                   \right.
\]
Note that we can rewrite \eqref{DPP} as the following deterministic form:
$$
 V\left(\hat{y}(\omega[n],\alpha[n])-\theta[n]\hat{T}(\omega[n],\alpha[n])\right)
 +\sum_{l=1}^LQ_l[n]\left(\hat{z}_l(\omega[n],\alpha[n])-c_l\hat{T}(\omega[n],\alpha[n])\right),
$$
Thus, Algorithm \ref{online-algorithm} proceeds by observing $\omega[n]$ on each frame $n$ and then choosing $\alpha[n]$ in $\mathcal{A}$ to minimize the above deterministic function.
We can now see that we only use knowledge of current realization $\omega[n]$, not statistics of $\omega[n]$. Also, the compactness assumption (Assumption \ref{assumption-for-algorithm}) guarantees that the minimum of \eqref{DPP} is always achievable.

\begin{algorithm}
\begin{itemize}
  \item At the beginning of each frame $n$, the controller observes $Q_l[n]$, $\theta[n]$, $\omega[n]$ and chooses action $\alpha[n]\in\mathcal{A}$ to minimize the following function:
\begin{equation}\label{DPP}
  \expect{\left. V(y[n]-\theta[n]T[n])+\sum_{l=1}^LQ_l[n](z_l[n]-c_lT[n])\right|Q_l[n],\theta[n],\omega[n]}.
\end{equation}
  \item Update $\theta[n]$:
  \[\theta[n+1]=\left[\frac{1}{(n+1)^{\delta}}\sum_{i=0}^{n}\left(y[i]-\theta[i]T[i]
  +\frac{1}{V}\sum_{l=1}^LQ_l[i](z_l[i]-c_lT[i])\right)\right]_{0}^{\theta_{\max}}.\]
  \item Update virtual queues $Q_l[n]$:
  \[Q_l[n+1]=\max\{Q_l[n]+z_l[n]-c_lT[n],0\},~l=1,2,\cdots,L.\]
\end{itemize}
\caption{Online renewal optimization:}
\label{online-algorithm}
\end{algorithm}

\section{Feasibility Analysis}
In this section, we prove that the proposed algorithm gives a sequence of actions $\{\alpha[n]\}_{n=0}^{\infty}$ which satisfies all desired constraints with probability 1. Specifically, we show that all virtual queues are stable with probability 1, in which we leverage an important lemma from \citet{hajek1982hitting} to obtain a exponential bound for the norm of $\mathbf{Q}[n]$.

\subsection{The drift-plus-penalty bound}
The start of our proof uses the drift-plus-penalty methodology. For a general introduction on this topic, see \citet{neely2012stability} for more details. 
We define the 2-norm function of the virtual queue vector as:
\[\|\mathbf{Q}[n]\|^2=\sum_{l=1}^LQ_l[n]^2.\]
Define the \textit{Lyapunov drift} $\Delta(\mathbf{Q}[n])$ as
\[\Delta(\mathbf{Q}[n])=\frac12\left(\|\mathbf{Q}[n+1]\|^2-\|\mathbf{Q}[n]\|^2\right).\]
Next, define the penalty function at frame $n$ as $V(y[n]-\theta[n]T[n])$, where $V>0$ is a fixed trade-off parameter. Then, the drift-plus-penalty methodology suggests that we can stabilize the virtual queues by choosing an action $\alpha[n]\in\mathcal{A}$ to greedily minimize the following drift-plus-penalty expression, with the observed $\mathbf{Q}[n]$, $\omega[n]$ and $\theta[n]$:
\[\expect{\left.V(y[n]-\theta[n]T[n])+\Delta(\mathbf{Q}[n])\right|Q_l[n],\theta[n],\omega[n]}.\]
The penalty term $V(y[n]-\theta[n]T[n])$ uses the $\theta[n]$ variable, which 
depends on events from all previous frames. This penalty does not
fit the rubric of \citet{neely2012stability} and convergence of the algorithm does
not follow from prior work. A significant thrust of the current paper is convergence
analysis under such a penalty function.

In order to obtain an upper bound on $\Delta(\mathbf{Q}[n])$, we square both sides of \eqref{queue-update} and use the fact that $\max\{x,0\}^2\leq x^2$,
\begin{align}\label{dpp-relation}
Q_l[n+1]^2\leq Q_l[n]^2+(z_l[n]-c_lT[n])^2+2Q_l[n](z_l[n]-c_lT[n]).
\end{align}
Summing the above over all $l \in \{1, \ldots, L\}$ and dividing by $2$ gives
$$ \Delta(\mathbf{Q}[n]) \leq \frac{1}{2}\sum_{l=1}^L (z_l[n]-c_lT[n])^2 + \sum_{l=1}^LQ_l[n](z_l[n]-c_lT[n])$$
Adding $V(y[n]-\theta[n]T[n])$ to both sides and taking conditional expectations gives
\begin{align}
&\expect{\left.V(y[n]-\theta[n]T[n])+\Delta(\mathbf{Q}[n])\right|Q_l[n],\theta[n],\omega[n]}\nonumber\\
\leq& \expect{\left.V(y[n]-\theta[n]T[n])+\sum_{l=1}^LQ_l[n](z_l[n]-c_lT[n])\right|Q_l[n],\theta[n],\omega[n]}
+\frac12\sum_{l=1}^L\expect{(z_l[n]-c_lT[n])^2}\nonumber\\
\leq& \expect{\left.V(y[n]-\theta[n]T[n])+\sum_{l=1}^LQ_l[n](z_l[n]-c_lT[n])\right|Q_l[n],\theta[n],\omega[n]}
+\frac{B^2}{\eta^2}.
\label{dpp-upperbound}
\end{align}
where the last inequality follows from Proposition \ref{prop-1}.
Thus, as we have already seen in Algorithm \ref{online-algorithm},
the proposed algorithm observes the vector $\mathbf{Q}[n]$, the random event $\omega[n]$ and $\theta[n]$ at frame $n$, and minimizes the right hand side of \eqref{dpp-upperbound}. 

\subsection{Bounds on the virtual queue process and feasibility}
In this section, we show how the bound \eqref{dpp-upperbound} leads to the feasibility of the proposed algorithm.
Define $\mathcal{H}_n$ as the system history information up until frame $n$. Formally, $\{\mathcal{H}_n\}_{n=0}^{\infty}$ is a filtration where each $\mathcal{H}_n$ is the $\sigma$-algebra generated by all the random variables before frame $n$. Notice that since $\mathbf{Q}[n]$ and $\theta[n]$ depend only on the events before frame $n$, $\mathcal{H}_n$ contains both $\mathbf{Q}[n]$ and $\theta[n]$.
 The following important lemma gives a stability criterion for any given real random process with certain negative drift property:

\begin{lemma}[Theorem 2.3 of \citet{hajek1982hitting}]\label{master-bound}
Let $R[n]$ be a real random process over $n\in \{0,1,2,\cdots\}$ satisfying the following two conditions for a fixed $r>0$:
\begin{enumerate}
\item For any $n$, $\expect{\left.e^{r(R[n+1]-R[n])}\right| \mathcal{H}_n}\leq \Gamma$, for some $\Gamma>0$.
\item Given $R[n]\geq\sigma$, $\expect{\left.e^{r(R[n+1]-R[n])}\right| \mathcal{H}_n}\leq \rho$, with some $\rho\in(0,1)$.
\end{enumerate}
Suppose further that $R[0]\in\mathbb{R}$ is given and finite, then, at every $n\in\{0,1,2,\cdots\}$, the following bound holds:
\[\expect{e^{rR[n]}}\leq \rho^ne^{rR[0]}+\frac{1-\rho^n}{1-\rho}\Gamma e^{r\sigma}.\]
\end{lemma}

Thus, in order to show the stability of the virtual queue process, it is enough to test the above two conditions with
$R[n]=\|\mathbf{Q}[n]\|$. The following lemma shows that $\|\mathbf{Q}[n]\|$ satisfies these two conditions:
\begin{lemma}[Drift condition]\label{geometric-bound}
Let $R[n]=\|\mathbf{Q}[n]\|$, then, it satisfies the two conditions in Lemma \ref{master-bound} with the following constants:
\begin{align*}
\Gamma&=B,\\
r&=\min\left\{\eta,\frac{\xi\eta^2}{4B}\right\},\\
\sigma&=C_0V,\\
\rho&=1-\frac{r\xi}{2}+\frac{2B}{\eta^2}r^2<1.
\end{align*}
where
$C_0=\frac{2B^2}{V\xi\eta^2}+\frac{2(\theta_{\max}+1)B}{\xi\eta}-\frac{\xi}{4V}$.
\end{lemma}
%
%

The central idea of the proof is to plug the $\xi$-slackness policy specified in Assumption \ref{slack} into the right hand side of \eqref{dpp-upperbound}. A similar idea has been presented in the
Lemma 6 of \citet{wei2015probabilistic} under the bounded increment of the virtual queue process. Here, we generalize the idea to the case where the increment of the virtual queues contains exponential type random variables $z_l[n]$ and $T[n]$. Note that
the boundedness of $\theta[n]$ is crucial for the argument to hold, which justifies the truncation of pseudo average in the algorithm. Lemma \ref{master-bound} is proved in the Appendix \ref{sec:proof}.

Combining the above two lemmas, we immediately have the following corollary:
\begin{corollary}[Exponential decay]\label{exponential-queue-bound}
Given $\mathbf{Q}[0]=0$, the following holds for any $n\in\{0,1,2,\cdots\}$ under the proposed algorithm,\\
\begin{equation}\label{eq:exp-queue-bound}
\expect{e^{r\|\mathbf{Q}[n]\|}}\leq D,
\end{equation}
where
\[D=1+\frac{B}{1-\rho} e^{rC_0V},\]
and $r,~\rho,~C_0$ are as defined in Lemma \ref{geometric-bound}. Furthermore, we have
$\expect{\|Q[n]\|}\leq \frac{1}{r}\log(1+ \frac{B}{1-\rho} e^{rC_0V})$, i.e. the queue size is $\mathcal{O}(V)$.
\end{corollary}
The bound on $\expect{\|Q[n]\|}$ follows readily from \eqref{eq:exp-queue-bound} via Jensen's inequality.
With Corollary \ref{exponential-queue-bound} in hand, we can prove the following theorem:
\begin{theorem}[Feasibility]\label{feasibility}
All constraints in \eqref{prob-1}-\eqref{prob-3} are satisfied under the proposed algorithm with probability 1.
\end{theorem}
\proof{Proof of Theorem \ref{feasibility}.~~}
By queue updating rule \eqref{queue-update}, for any $n$ and any $l\in\{1,2,\cdots,L\}$, one has
  \[Q_l[n+1]\geq Q_l[n]+z_l[n]-c_lT[n].\]
Fix $N$ as a positive integer. Then, summing over all $n\in\{0,1,2,\cdots,N-1\}$,
  \[Q_l[N]\geq Q_l[0]+\sum_{n=0}^{N-1}(z_l[n]-c_lT[n]).\]
  Since $Q_l[0]=0,~\forall l$ and $T[n]\geq1,~\forall n$,
  \begin{equation}\label{inter-constraint-violation}
  \frac{\sum_{n=0}^{N-1}z_l[n]}{\sum_{n=0}^{N-1}T[n]}-c_l\leq\frac{Q_l[N]}{\sum_{n=0}^{N-1}T[n]}\leq\frac{Q_l[N]}{N}.
  \end{equation}
Define the event
\[A_{N}^{(\varepsilon)}=\{Q_l[N]>\varepsilon N\}.\]
By the Markov inequality and Corollary \ref{exponential-queue-bound}, for any $\varepsilon>0$, we have
\begin{align*}
Pr(Q_l[N]>\varepsilon N)
\leq&Pr\left(r\|\mathbf Q[N]\|> r\varepsilon N\right)\\
=&Pr\left(e^{r\|\mathbf Q[N]\|}> e^{r\varepsilon N}\right)\\
\leq&\frac{\expect{e^{r\|\mathbf Q[N]\|}}}{e^{r\varepsilon N}}\leq De^{-r\varepsilon N},
\end{align*}
where $r$ is defined in Corollary \ref{exponential-queue-bound}.
Thus, we have
\begin{align*}
\sum_{N=0}^{\infty}Pr(Q_l[N]>\varepsilon N)\leq D\sum_{N=0}^{\infty}e^{-r\varepsilon N}<+\infty.
\end{align*}
Thus, by the Borel-Cantelli lemma \citet{Durrett},
\[Pr\left(A_{N}^{(\varepsilon)}~\textrm{occurs infinitely often}\right)=0.\]
Since $\varepsilon>0$ is arbitrary, letting $\varepsilon\rightarrow0$ gives
\[Pr\left(\lim_{N\rightarrow\infty}\frac{Q_l[N]}{N}=0\right)=1.\]
Finally, taking the $\limsup_{N\rightarrow\infty}$ from both sides of \eqref{inter-constraint-violation} and substituting in the  above equation gives the claim.
\Halmos
\endproof

\section{Optimality Analysis}
In this section, we show that the proposed algorithm achieves time average penalty within $\mathcal{O}(1/V)$ of the optimal objective $\theta^*$. Since the algorithm meets all the constraints, it follows,
\[\limsup_{n\rightarrow\infty}\frac{\sum_{i=0}^{n-1}y[i]}{\sum_{i=0}^{n-1}T[i]}\geq\theta^*,~~w.p.1.\]
Thus, it is enough to prove the following theorem:
\begin{theorem}[Near optimality]\label{theorem_average_converge}
For any $\delta\in(1/3,1)$ and $V\geq1$, the objective value produced by the proposed algorithm is near optimal with
\[\limsup_{n\rightarrow\infty}\frac{\sum_{i=0}^{n-1}y[i]}{\sum_{i=0}^{n-1}T[i]}\leq\theta^*+\frac{B^2}{\eta^2V},~w.p.1,\]
i.e. the algorithm achieves $\mathcal{O}(1/V)$ near optimality.
\end{theorem}
\begin{remark}
Combining Theorem \ref{theorem_average_converge} with Corollary \ref{exponential-queue-bound}, we see that the tuning parameter $V$ plays a trade-off between the sub-optimality and the virtual queue bound (i.e. the constraint violation). In particular, our result recovers the classical 
$[\mathcal{O}(1/V),~\mathcal{O}(V)]$ trade-off in the work of opportunistic scheduling \citet{Neely2010}. 
\end{remark}

In order to prove Theorem \ref{theorem_average_converge}, we introduce the following notation:
\begin{align*}
&\textrm{original pseudo average}:~~\hat{\theta}[n]\triangleq\frac{1}{(n+1)^\delta}\sum_{i=0}^{n}\left(y[i]-\theta[i]T[i]+\frac{1}{V}\sum_{l=1}^LQ_l[i](z_l[i]-c_lT[i])\right),\\
&\textrm{tamed pseudo average}:~~\theta[n]\triangleq\left[\frac{1}{(n+1)^\delta}\sum_{i=0}^{n}\left(y[i]-\theta[i]T[i]+\frac{1}{V}\sum_{l=1}^LQ_l[i](z_l[i]-c_lT[i])\right)\right]_0^{\theta_{\max}}.
\end{align*}

\subsection{Relation between $\hat\theta[n]$ and $\theta[n]$}
We start with a preliminary lemma illustrating that the original pseudo average $\hat\theta[n]$ behaves almost the same as the tamed pseudo average $\theta[n]$. Note that $\theta[n]$ can be written as:
$$\theta[n] = [ \hat{\theta}[n]]_0^{\theta_{max}}. $$
 
\begin{lemma}[Equivalence relation]\label{properties}
For any $x\in(0,\theta_{\max})$,
\begin{enumerate}
\item $\theta[n]\geq x$ if and only if $\hat{\theta}[n]\geq x$.
\item $\theta[n]\leq x$ if and only if $\hat{\theta}[n]\leq x$.
\item $\limsup_{n\rightarrow\infty}\theta[n]\leq x$ if and only if $\limsup_{n\rightarrow\infty}\hat\theta[n]\leq x$.
\item $\limsup_{n\rightarrow\infty}\theta[n]\geq x$ if and only if $\limsup_{n\rightarrow\infty}\hat\theta[n]\geq x$.
\end{enumerate}
\end{lemma}

This lemma is intuitive and the proof is shown in the Appendix \ref{sec:proof}. We will prove results on $\hat{\theta}[n]$ which extend naturally to $\theta[n]$ via Lemma \ref{properties}.

The key idea of proving Theorem \ref{theorem_average_converge} is to bound the original pseudo average process $\hat{\theta}[n]$ asymptotically from above by $\theta^*$, which is Theorem \ref{theorem-asymptotic-upperbound} below. We then prove Theorem \ref{theorem-asymptotic-upperbound} through the following three steps:
\begin{itemize}
\item We construct a truncated version of $\hat\theta[n]$, namely $\tilde\theta[n]$, which has the same limit as $\hat{\theta}[n]$ (Lemma \ref{truncation-lemma} below), so that it is enough to show $\tilde\theta[n]\leq \theta^*$ asymptotically.
\item  For the process $\tilde\theta[n]$, we bound the moments of the hitting time, namely, the time interval between two consecutive visits to the region $\{\tilde\theta[n]\leq\theta^*\}$, by constructing a dominating exponential supermartingale and bounding its size.
 (Lemma \ref{exp-supMG} and \ref{time-moment-bound} below).

\item
We show that $\tilde\theta[n]>\theta^*$ only finitely often asymptotically (with probability 1) using the bounded moments of the hitting time.  
\end{itemize}

\subsection{Towards near optimality (I): Truncation}
The following lemma states that  the optimality of \eqref{prob-1}-\eqref{prob-3} is achievable within the closure of the set of all one-shot averages specified in Assumption \ref{compact}:
\begin{lemma}[Stationary optimality]  \label{optimal-stationary-lemma}
Let $\theta^*$ be the optimal objective of \eqref{prob-1}-\eqref{prob-3}. Then, there exists a tuple $(y^*,~T^*,~\mathbf{z}^*)\in\overline{\mathcal{R}}$, the closure of $\mathcal{R}$, such that the following hold:
\begin{align}
&y^*/T^*=\theta^*   \label{iid1}\\
&z_l^*/T^*\leq c_l ,~\forall l\in\{1,2,\cdots,L\}, \label{iid2}
\end{align}
i.e. the optimality is achievable within $\overline{\mathcal{R}}$.
\end{lemma}
The proof of this lemma is similar to the proof of Theorem 4.5 as well as Lemma 7.1 of \citet{Neely2010}. We omit the details for brevity.

We start the truncation by picking up an $\varepsilon_0>0$ small enough so that $\theta^*+\varepsilon_0/V<\theta_{\max}$.
We aim to show $\limsup_{n\rightarrow\infty}\theta[n]\leq\theta^*+\varepsilon_0/V$. By Lemma \ref{properties},
it is enough to show $\limsup_{n\rightarrow\infty}\hat{\theta}[n]\leq\theta^*+\varepsilon_0/V$. The following lemma tells us it is enough to prove it on a further term-wise truncated version of $\hat{\theta}[n]$.

\begin{lemma}[Truncation lemma]\label{truncation-lemma}
Consider the following alternative pseudo average $\{\tilde{\theta}[n]\}_{n=0}^{\infty}$ obtained by truncating each 
summand such that $\tilde{\theta}[0]=0$ and
\begin{align*}
\tilde{\theta}[n+1]=\frac{1}{(n+1)^{\delta}}\sum_{i=0}^n\left[\left(y[i]-\theta[i]T[i]+\frac1V\sum_{l=1}^LQ_l[i](z_l[i]-c_lT[i])\right)
\wedge \left(\left(\frac{2}{\eta}+\frac{4\sqrt{L}}{\eta rV}\right)\log^2(i+1)\right)\right],
\end{align*}
where $a\wedge b\triangleq\min\{a,b\}$, $\eta$ is defined in Assumption \ref{bounded-assumption}
and $r$ is defined in Lemma \ref{geometric-bound}. Then, we have
\[\limsup_{n\rightarrow\infty}\hat{\theta}[n]=\limsup_{n\rightarrow\infty}\tilde{\theta}[n].\]
\end{lemma}
\proof{Proof of Lemma \ref{truncation-lemma}.~~}
Consider any frame $i \in \{0, 1, 2, \ldots\}$ such that there is a discrepancy between the summand of $\hat{\theta}[n]$ and $\tilde{\theta}[n]$, i.e.
\begin{align}\label{discrepancy-condition}
y[i]-\theta[i]T[i]+\frac1V\sum_{l=1}^LQ_l[i](z_l[i]-c_lT[i])>\left(\frac{2}{\eta}+\frac{4\sqrt{L}}{\eta rV}\right)\log^2(i+1),
\end{align}
By the Cauchy-Schwartz inequality, this implies
\begin{align*}
y[i]-\theta[i]T[i]+\frac1V\sqrt{\sum_{l=1}^LQ_l[i]^2}\sqrt{\sum_{l=1}^L(z_l[i]-c_lT[i])^2}>\left(\frac{2}{\eta}+\frac{4\sqrt{L}}{\eta rV}\right)\log^2(i+1).
\end{align*}
Thus, at least one of the following three events happened:
\begin{enumerate}
\item $A_i\triangleq \left\{y[i]-\theta[i]T[i]>\frac{2}{\eta}\log^2(i+1)\right\}$.
\item $B_i\triangleq\left\{\sqrt{\sum_{l=1}^LQ_l[i]^2}>\frac{2\sqrt{L}}{r}\log(i+1)\right\}$.
\item $E_i\triangleq\left\{K[i]>\frac{2}{\eta}\log(i+1)\right\}$.
\end{enumerate}
where $K[i]$ is defined in \eqref{def-K}. Indeed, the occurence of one of the three events is necessary for \eqref{discrepancy-condition} to happen.
We then argue that these three events jointly occur only finitely many times. Thus, as $n\rightarrow\infty$, the discrepancies are negligible.

Assume the event $A_i$ occurs, then, since $y[i]-\theta[i]T[i]\leq y[i]$, it follows
$y[i]>\frac{2}{\eta}\log^2(i+1)$.
Then, we have
\begin{align*}
Pr(A_i)\leq&Pr\left(y[i]>\frac{2}{\eta}\log^2(i+1)\right)\\
=&Pr\left(e^{\eta y[i]}>e^{2\log^2(i+1)}\right)\\
\leq&\frac{\expect{e^{\eta y[i]}}}{(i+1)^{2\log(i+1)}}\leq\frac{B}{(i+1)^{2\log(i+1)}},
\end{align*}
where the second to last inequality follows from the Markov inequality  
and the last inequality follows from Assumption \ref{bounded-assumption}.

Assume the event $B_i$ occurs, then, we have
\begin{align*}
\|\mathbf{Q}[i]\|=\sqrt{\sum_{l=1}^LQ_l[i]^2}>\frac{2\sqrt{L}}{r}\log(i+1)\geq\frac{2}{r}\log(i+1).
\end{align*}
Thus, 
\begin{align*}
Pr(B_i)\leq& Pr\left(\|\mathbf{Q}[i]\|>\frac{2}{r}\log(i+1)\right)\\
=&Pr\left(e^{r\|\mathbf{Q}[i]\|}>e^{2\log (i+1)}\right)\\
\leq&\frac{\expect{e^{r\|\mathbf{Q}[i]\|}}}{(i+1)^2}\leq\frac{D}{(i+1)^2},
\end{align*}
where the second to last inequality follows from the 
Markov inequality and the last inequality follows from Corollary \ref{exponential-queue-bound}.

Assume the event $E_i$ occurs.
Again, by Assumption \ref{bounded-assumption} and the Markov inequality,
\begin{align*}
Pr(E_i)=&Pr\left(K[i]>\frac{2}{\eta}\log(i+1)\right)\\
=&Pr\left(e^{\eta K[i]}>e^{2\log(i+1)}\right)\\
\leq&\frac{\expect{e^{\eta K[i]}}}{(i+1)^2}\leq\frac{B}{(i+1)^2},
\end{align*}
where the last inequality follows from Assumption \ref{bounded-assumption} again.
Now, by a union bound,
$$Pr(A_i\cup B_i\cup E_i)\leq Pr(A_i)+Pr(B_i)+Pr(E_i)\leq\frac{B}{(i+1)^{2\log(i+1)}}+\frac{B+D}{(i+1)^2},$$
and thus,
$$\sum_{i=0}^{\infty}Pr(A_i\cup B_i\cup E_i)\leq\sum_{i=0}^{\infty}
\left(\frac{B}{(i+1)^{2\log(i+1)}}+\frac{B+D}{(i+1)^2}\right)<\infty$$
By the Borel-Cantelli lemma, we have the joint event $A_i\cup B_i\cup E_i$ occurs only finitely many times with probability 1, and our proof is finished.
\Halmos
\endproof

Lemma \ref{truncation-lemma} is crucial for the rest of the proof. Specifically, it creates an alternative sequence $\tilde{\theta}[n]$ which has the following two properties:
\begin{enumerate}
\item We know exactly what the upper bound of each of the summands is, whereas in $\hat{\theta}[n]$, there is no exact bound for the summand due to $Q_l[i]$ and other exponential type random variables.
\item For any $n\in\mathbb{N}$, we have $\tilde{\theta}[n]\leq\hat{\theta}[n]$. Thus, if $\tilde{\theta}[n]\geq\theta^*+\varepsilon_0/V$ for some $n$, then, $\hat{\theta}[n]\geq\theta^*+\varepsilon_0/V$.
\end{enumerate}

\subsection{Towards near optimality (II): Exponential supermartingale}
The following preliminary lemma demonstrates a negative drift property for each of the summands in $\tilde{\theta}[n]$.
\begin{lemma}[Key feature inequality]\label{key-feature}
For any $\varepsilon_0>0$, if $\theta[i]\geq\theta^*+\varepsilon_0/V$, then, we have
\begin{align*}
\expect{\left.\left(y[i]-\theta[i]T[i]+\frac1V\sum_{l=1}^LQ_l[i](z_l[i]-c_lT[i])\right)
\wedge\left(\left(\frac{2}{\eta}+\frac{4\sqrt{L}}{\eta rV}\right)\log^2(i+1)\right)\right|\mathcal{H}_i}\leq-\varepsilon_0/V,
\end{align*}
\end{lemma}
\proof{Proof of Lemma \ref{key-feature}.~~}
Since the proposed algorithm minimizes \eqref{DPP} over all possible decisions in $\mathcal{A}$, it must achieve value less than or equal to that of any randomized stationary algorithm $\alpha^*[i]$. This in turn implies,
\begin{align*}
&\expect{\left.\left(y[i]-\theta[i]T[i]+\frac1V\sum_{l=1}^LQ_l[i](z_l[i]-c_lT[i])\right)\right|\mathcal{H}_i,\omega[i]}\\
\leq&\expect{\left.\left(\hat{y}(\omega[i],\alpha^*[i])-\theta[i]\hat T(\omega[i],\alpha^*[i])
+\frac1V\sum_{l=1}^LQ_l[i](\hat z_l(\omega[i],\alpha^*[i])-c_l\hat T(\omega[i],\alpha^*[i]))\right)\right|\mathcal{H}_i,\omega[i]}.
\end{align*}
Taking expectation from both sides with respect to $\omega[i]$ and using the fact that randomized stationary algorithms are i.i.d. over frames and independent of $\mathcal{H}_i$, we have
\begin{align*}
&\expect{\left.\left(y[i]-\theta[i]T[i]+\frac1V\sum_{l=1}^LQ_l[i](z_l[i]-c_lT[i])\right)\right|\mathcal{H}_i,\omega[i]}
\leq\overline{y}-\theta[i]\overline T
+\frac1V\sum_{l=1}^LQ_l[i](\overline z_l-c_l\overline T),
\end{align*}
for any $(\overline{y},\overline{T},\overline{\mathbf{z}})\in\mathcal{R}$. Since $(y^*,T^*,\mathbf{z}^*)$  specified in Lemma \ref{optimal-stationary-lemma} is in the closure of $\mathcal{R}$, we can replace $(\overline{y},\overline{T},\overline{\mathbf{z}})$ by the tuple $(y^*,T^*,\mathbf{z}^*)$ and the inequality still holds. This gives
\begin{align*}
&\expect{\left.\left(y[i]-\theta[i]T[i]+\frac1V\sum_{l=1}^LQ_l[i](z_l[i]-c_lT[i])\right)\right|\mathcal{H}_i,\omega[i]}\\
\leq& y^*-\theta[i]T^*+\frac1V\sum_{l=1}^LQ_l[i](z_l^*-c_l T^*),\\
=& T^*\left(y^*/T^*-\theta[i]+\frac1V\sum_{l=1}^LQ_l[i](z_l^*/T^*-c_l)\right)\\
\leq&T^*(\theta^*-\theta[i])\leq-\varepsilon_0/V,
\end{align*}
where the second to last inequality follows from \eqref{iid1} and \eqref{iid2}, and the last inequality follows from $\theta[i]\geq\theta^*+\varepsilon_0/V$ and $T[i]\geq1$.
Finally, since $a\wedge b\leq a$ for any real numbers $a,b$, it follows,
\begin{align*}
&\expect{\left.\left(y[i]-\theta[i]T[i]+\frac1V\sum_{l=1}^LQ_l[i](z_l[i]-c_lT[i])\right)
\wedge\left(\left(\frac{2}{\eta}+\frac{4\sqrt{L}}{\eta rV}\right)\log^2(i+1)\right)\right|\mathcal{H}_i}\\
\leq&
\expect{\left.\left(y[i]-\theta[i]T[i]+\frac1V\sum_{l=1}^LQ_l[i](z_l[i]-c_lT[i])\right)\right|\mathcal{H}_i}
\leq-\varepsilon_0/V,
\end{align*}
and the claim follows.
\Halmos
\endproof

Define $n_k$ as the frame where $\tilde{\theta}[n]$ visits the set $(-\infty,~\theta^*+\varepsilon_0/V)$ for the $k$-th time with the following conventions: 1. If $\tilde{\theta}[n]\in(-\infty,~\theta^*+\varepsilon_0/V)$ and $\tilde{\theta}[n+1]\in(-\infty,~\theta^*+\varepsilon_0/V)$, then we count them as two times. 2. When $k=1$, $n_1$ is equal to 0. Define the \textit{hitting time} $S_{n_k}$ as
\[S_{n_k}=n_{k+1}-n_k.\]
The goal is to obtain a moment bound on this quantity when $\tilde\theta[n_k+1]\geq\theta^*+\varepsilon_0/V$ (otherwise, this quantity is 1). In order to do so, we introduce a new process as follows. For any $n_k$, define 
\begin{align}\label{Fn-construction}
F[n]\triangleq\sum_{i=n_k}^{n-1}\left(y[i]-\theta[i]T[i]+\frac1V\sum_{l=1}^LQ_l[i](z_l[i]-c_lT[i])\right)
\wedge\left(\left(\frac{2}{\eta}+\frac{4\sqrt{L}}{\eta rV}\right)\log^2(i+1)\right),~\forall n>n_k,
\end{align}

The following lemma shows that indeed this $F[n]$ is closely related to $\tilde\theta[n]$. It plays an important role in proving Lemma \ref{time-moment-bound}:
\begin{lemma}\label{comparison-lemma}
~~For any $n>n_k$, if $\tilde{\theta}[n]\geq\theta^*+\varepsilon_0/V$, then, $F[n]\geq0$.
\end{lemma}
\proof{Proof of Lemma \ref{comparison-lemma}.~~}
Suppose $\tilde{\theta}[n]\geq\theta^*+\varepsilon_0/V$, then, the following holds
\[\theta^*+\varepsilon_0/V\leq\tilde{\theta}[n]=\frac{n_k^\delta}{n^\delta}\tilde\theta[n_k]+\frac{1}{n^\delta}F[n].\]
Thus,
$$F[n]\geq n^\delta(\theta^*+\varepsilon_0/V)-n_k^\delta\tilde\theta[n_k].$$
Since at the frame $n_k$, $\tilde\theta[n_k]<\theta^*+\varepsilon_0/V$, it follows,
\[F[n]\geq \left(n^\delta-n_k^\delta\right)(\theta^*+\varepsilon_0/V).\]
Since $\theta^*+\varepsilon_0/V\geq0$, it follows $F[n]\geq0$ and the claim follows.
\Halmos
\endproof
%
%
 Recall our goal is to bound the hitting time $S_{n_k}$ of the process $\tilde\theta[n]$ when $\{\tilde\theta[n_k+1]\geq\theta^*+\varepsilon_0/V\}$, with a strictly negative drift property as Lemma \ref{key-feature}. A classical approach analyzing the hitting time of a stochastic process came from Wald's construction of martingale for sequential analysis (see, for example,  \citet{wald1944cumulative} for details). Later,
  \citet{hajek1982hitting} extended this idea to analyze the stability of a queueing system with drift condition by a supermartingale construnction. Here, we take one step further by considering the following supermartingale construction based on $F[n]$:
\begin{lemma}[Exponential Supermartingale]\label{exp-supMG}
Fix $\varepsilon_0>0$ and $V\geq\max\left\{\frac{\varepsilon_0\eta}{4\log^22}-\frac{2\sqrt{L}}{r},~1\right\}$ such that $\theta^*+\varepsilon_0/V<\theta_{\max}$.
Define a new random process $G[n]$ starting from $n_k+1$ with
\[G[n]\triangleq\frac{\exp\left(\lambda_nF[n\wedge(n_k+S_{n_k})]\right)}{\prod_{i=n_k+1}^{n\wedge (n_k+S_{n_k})}\rho_i}\mathbf{1}_{\{\tilde\theta[n_k+1]\geq\theta^*+\varepsilon_0/V\}},\]
where for any set $A$, $\mathbf{1}_A$ is the indicator function which takes value 1 if $A$ is true and 0 otherwise.
For any $n\geq n_k+1$, $\lambda_n$ and $\rho_n$ are defined as follows:
\begin{align*}
\lambda_n=&\frac{\varepsilon_0}{2Ve\left(\frac{2}{\eta}+\frac{4\sqrt{L}}{\eta rV}\right)^2\log^4(n+1)},\\
\rho_n=&1-\frac{\varepsilon_0^2}{4V^2e\left(\frac{2}{\eta}+\frac{4\sqrt{L}}{\eta rV}\right)^2\log^4(n+1)}.
\end{align*}
Then, the process $G[n]$ is measurable with respect to $\mathcal{H}_n$, $\forall n\geq n_k+1$, and furthermore, it is a supermartingale with respect to the filtration $\{\mathcal{H}_n\}_{n\geq n_k+1}$.
\end{lemma}
The proof of Lemma \ref{exp-supMG} is shown in Appendix \ref{sec:proof}.
\begin{remark}
~~If the increments $F[n+1]-F[n]$ were to be bounded, then, we could adopt the similar construction as that of \citet{hajek1982hitting}. However, in our scenario $F[n+1]-F[n]$ is of the order $\log^2(n+1)$, which is increasing and unbounded. Thus, we need decreasing exponents $\lambda_n$ and increasing weights $\rho_n$ to account for that. Furthermore, the indicator function indicates that we are only interested in the scenario $\{\tilde\theta[n_k+1]\geq\theta^*+\varepsilon_0/V\}$.
\end{remark}

The following lemma uses the previous result to bound the conditional fourth moment of the hitting time $S_{n_k}$. 
\begin{lemma}\label{time-moment-bound}
~~Given
$V\geq\max\left\{\frac{\varepsilon_0\eta}{4\log^22}-\frac{2\sqrt{L}}{r},~1\right\}$ as in Lemma \ref{exp-supMG}, for any $\beta\in(0,1/5)$ and any $\varepsilon_0>0$ such that $\theta^*+\varepsilon_0/V<
\theta_{\max}$,
there exists a positive constant $C_{\beta,V,\varepsilon_0}\simeq\mathcal{O}\left(V^{10}\beta^{-20}\varepsilon_0^{-10}\right)$, such that
\[\expect{S_{n_k}^4|\mathcal{H}_{n_k}}\leq C_{\beta,V,\varepsilon_0}(n_k+2)^{4\beta},~~\forall k\geq1.\]
\end{lemma}

\proof{Proof of Lemma \ref{time-moment-bound}.~~}
First of all, from Lemma \ref{exp-supMG} gives that $G[n]$ is a supermartingale starting from $n_k+1$, thus, we have the following chains of inequalities for any $n\geq n_k+1$:
\begin{align*}
G[n_k+1]=&\expect{G[n_k+1]~|~\mathcal{H}_{n_k+1}}\\
\geq&\expect{G[n]~|~\mathcal{H}_{n_k+1}}\\
=&\expect{\left.\frac{e^{\lambda_n F[n\wedge(n_k+S_{n_k})]}}
{\prod_{i=n_k+1}^{n}\rho_i}\mathbf{1}_{\{\tilde\theta[n_k+1]\geq\theta^*+\varepsilon_0/V\}}~\right|~\mathcal{H}_{n_k+1}}\\
\geq&\expect{\left.\frac{e^{\lambda_n F[n\wedge(n_k+S_{n_k})]}}
{\prod_{i=n_k+1}^{n}\rho_i}\mathbf{1}_{\{S_{n_k}\geq n-n_k+1\}}\mathbf{1}_{\{\tilde\theta[n_k+1]\geq\theta^*+\varepsilon_0/V\}}~\right|~\mathcal{H}_{n_k+1}}\\
\geq&\frac{1}
{\prod_{i=n_k+1}^{n}\rho_i}Pr\left[\left.S_{n_k}\geq n-n_k+1,~\tilde\theta[n_k+1]\geq\theta^*+\varepsilon_0/V~\right|~\mathcal{H}_{n_k+1}\right],
\end{align*}
where the first inequality uses the supermartingale property and the
 last inequality uses Lemma \ref{comparison-lemma} that on the set $\{S_{n_k}\geq n-n_k+1\}$, $n\wedge(n_k+S_{n_k})=n$ and $F[n]\geq0$. By definition of $G[n_k+1]$,
\begin{align*}
G[n_k+1]=\frac{e^{\lambda_{n_k+1} F[n_k+1]}}
{\rho_{n_k+1}}
\leq\frac{e^{\lambda_{n_k+1}\left(\frac{2}{\eta}+\frac{4\sqrt{L}}{\eta rV}\right)\log^2(n_k+2)}}{\rho_{n_k+1}}
\leq\frac43e,
\end{align*}
where the first inequality follows from the definition of $F[n]$, and the second inequality follows from the assumption that $V\geq\frac{\varepsilon_0\eta}{4\log^22}-\frac{2\sqrt{L}}{r}$, thus,
$\lambda_{n_k+1}\leq\frac{1}{\left(\frac{2}{\eta}+\frac{4\sqrt{L}}{\eta rV}\right)\log^2(n_k+2)}$ and $\rho_{n_k+1}\geq1-\frac{\log^22}{2e}>\frac34$. Thus,
it follows,
\[Pr\left[\left.S_{n_k}\geq n-n_k+1,~\tilde\theta[n_k+1]\geq\theta^*+\varepsilon_0/V~\right|~\mathcal{H}_{n_k+1}\right]\leq \left(\prod_{i=n_k+1}^{n}\rho_i\right)\cdot \frac{4}{3}e.\]
Now, we bound the fourth moment of hitting time:
\begin{align*}
&\expect{\left. S_{n_k}^4~\right|~\mathcal{H}_{n_k+1}}\\
=&\sum_{m=1}^{\infty}m^4Pr\left[\left.S_{n_k}= m~\right|~\mathcal{H}_{n_k+1}\right]\\
\leq&\sum_{m=1}^{\infty}\left((m+1)^4-m^4\right)Pr\left[\left.S_{n_k}\geq m+1,~\tilde\theta[n_k+1]\geq\theta^*+\varepsilon_0/V~\right|~\mathcal{H}_{n_k+1}\right]+1\\
\leq&4\sum_{m=1}^{\infty}(m+1)^3Pr\left[\left.S_{n_k}\geq m+1,~\tilde\theta[n_k+1]\geq\theta^*+\varepsilon_0/V~\right|~\mathcal{H}_{n_k+1}\right]+1\\
\leq&1+\frac{16}{3}e\sum_{m=1}^\infty(m+1)^3\prod_{i=n_k+1}^{n_k+m}\rho_i.
\end{align*}
Thus, it remains to
show there exists a constant $C$ on the order $\mathcal{O}\left(V^{10}\beta^{-20}\varepsilon_0^{-10}\right)$ such that
\[\sum_{m=1}^\infty(m+1)^3\prod_{i=n_k+1}^{n_k+m}\rho_i\leq C(n_k+2)^{4\beta},\]
which is given is Appendix \ref{computation}. This implies there exists a $C_{\beta,V,\varepsilon_0}$ so that
\[\expect{\left. S_{n_k}^4\right|\mathcal{H}_{n_k+1}}\leq C_{\beta, V,\varepsilon_0}(n_k+2)^{4\beta}.\]
Thus,
\begin{align*}
\expect{\left. S_{n_k}^4\right|\mathcal{H}_{n_k}}&=\expect{\expect{\left. S_{n_k}^4\right|\mathcal{H}_{n_k+1}}|\mathcal{H}_{n_k}}
\leq\expect{C_{\beta, V,\varepsilon_0}(n_k+2)^{4\beta}|\mathcal{H}_{n_k}}
=C_{\beta, V,\varepsilon_0}(n_k+2)^{4\beta},
\end{align*}
where the last equality follows from the fact that $n_k\in\mathcal{H}_{n_k}$.
This finishes the proof. 
\Halmos
\endproof

\subsection{An asymptotic upper bound on $\theta[n]$}\label{sec:finish-proof}
So far, we have proved that if we pick any $\varepsilon_0>0$ such that $\theta^*+\varepsilon_0/V<\theta_{\max}$, then,  the inter-visiting time has bounded conditional fourth moment. We aim to show that $\limsup_{n\rightarrow\infty}\hat\theta[n]\leq\theta^*$ with probability 1. By Lemma \ref{truncation-lemma}, it is enough to show $\limsup_{n\rightarrow\infty}\tilde\theta[n]\leq\theta^*$. To do so, we need the following Second Borel-Cantelli lemma:
\begin{lemma}[Theorem 5.3.2. of \citet{Durrett}]\label{second-borel}
Let $\mathcal{F}_k,~k\geq1$ be a filtration with $\mathcal{F}_1=\{\emptyset,\Omega\}$, and $A_k,~k\geq1$ be a sequence of events with $A_k\in\mathcal{F}_{k+1}$, then
\[\{A_k~\textrm{occurs infinitely often}\}=\left\{\sum_{k=1}^{\infty}Pr(A_k|\mathcal{F}_{k})=\infty\right\}\]
\end{lemma}

\begin{theorem}[Asymptotic upper bound]\label{theorem-asymptotic-upperbound}
For any 
$\delta\in(1/3,1)$ and $V\geq1$, the following hold,
\[\limsup_{n\rightarrow\infty}\hat\theta[n]\leq\theta^*,~~w.p.1,\]
and
\[\limsup_{n\rightarrow\infty}\theta[n]\leq\theta^*,~~w.p.1.\]
\end{theorem}
\proof{Proof of Theorem \ref{theorem-asymptotic-upperbound}.~~}
First of all, since the inter-hitting time $S_{n_k}$ has finite fourth moment, each inter-hitting time is finite with probability 1, and thus the process $\{\tilde{\theta}[n]\}_{n=0}^{\infty}$ will visit $(-\infty,\theta^*+\varepsilon_0/V)$ infinitely many times with probability 1.
Then, we pick any $\epsilon>0$ and define the following sequence of events:
\begin{equation}\label{def-A}
A_k\triangleq\left\{\frac{S_{n_k}}{n_k^{1/3}}>\epsilon\right\},~k=1,2,\cdots.
\end{equation}
For any fixed $k$, by Conditional Markov inequality, the following holds with probability 1:
\begin{align*}
Pr(A_k|\mathcal{H}_{n_k})=&Pr\left(\left.S_{n_k}^4>\epsilon^4 n_k^{4/3}\right|\mathcal{H}_{n_k}\right)\\
\leq&\frac{\expect{S_{n_k}^4|\mathcal{H}_{n_k}}}{\epsilon^4 n_k^{4/3}}\\
\leq&\frac{C_{\beta,V,\varepsilon_0}(n_k+2)^{4\beta}}{\epsilon^4 n_k^{4/3}}\\
\leq&\frac{C_{\beta,V,\varepsilon_0}}{\epsilon^4}n_k^{-4/3+4\beta}+\frac{C_{\beta,V,\varepsilon_0}2^{4\beta}}{\epsilon^4n_k^{4/3}}\\
\leq&\frac{C_{\beta,V,\varepsilon_0}}{\epsilon^4}k^{-4/3+4\beta}+\frac{C_{\beta,V,\varepsilon_0}2^{4\beta}}{\epsilon^4}k^{-4/3},
\end{align*}
where the second inequality follows from Lemma \ref{time-moment-bound} with $\beta\in(0,1/5)$, the third inequality follows from the fact that $(a+b)^x\leq a^x+b^x,~\forall a,b\geq0$ and $x\in(0,1)$. The last inequality follows from the fact that the inter-hitting time takes at least one frame and thus $n_k\geq k$.

Choose $\mathcal{F}_k=\mathcal{H}_{n_k}$ and $A_k$ as is defined in \eqref{def-A}. Then, for any $\beta\in(0,1/12)$, we have with probability 1,
\begin{align*}
\sum_{k=1}^{\infty}Pr(A_k|\mathcal{H}_{n_k})\leq
\sum_{k=1}^{\infty}\left(\frac{C_{\beta,V,\varepsilon_0}}{\epsilon^4}k^{-4/3+4\beta}+\frac{C_{\beta,V,\varepsilon_0}2^{4\beta}}{\epsilon^4}k^{-4/3}
\right)<\infty.
\end{align*}
Now by Lemma \ref{second-borel},
\[Pr\left(A_k~\textrm{occurs infinitely often}\right)=0.\]
Since the process $\{\tilde{\theta}[n]\}_{n=0}^{\infty}$ visits $(-\infty,\theta^*+\varepsilon_0/V)$ infinitely many times with probability 1, 
\[\limsup_{n\rightarrow\infty}\frac{S_{n_k}}{n_k^{1/3}}
=\limsup_{k\rightarrow\infty}\frac{S_{n_k}}{n_k^{1/3}}\leq\epsilon,~w.p.1,\]
Since $\epsilon>0$ is arbitrary, let $\epsilon\rightarrow0$ gives
\begin{equation}\label{bound-on-returning}
\lim_{n\rightarrow\infty}\frac{S_{n_k}}{n_k^{1/3}}=0,~w.p.1.
\end{equation}
Finally, we show how this convergence result leads to the bound of $\tilde{\theta}[n]$. According to the updating rule of $\tilde{\theta}[n]$, for any frame $n$ such that $n_k<n\leq n_{k+1}$,
\begin{align*}
\tilde{\theta}[n]=&(\frac{n_k}{n})^{\delta}\tilde{\theta}[n_k]
               +\frac{1}{n^{\delta}}\sum_{i=n_k}^{n-1}\left(y[i]-\theta[i]T[i]+\frac{1}{V}Q[i](z[i]-cT[i])
      \right)\wedge\left(\left(\frac{2}{\eta}+\frac{4\sqrt{L}}{\eta rV}\right)\log^2(i+1)\right)\\
\leq&(\frac{n_k}{n})^{\delta}\left(\theta^*+\frac{\varepsilon_0}{V}\right)
      +\frac{1}{n^{\delta}}\sum_{i=n_k}^{n-1}\left(\left(\frac{2}{\eta}+\frac{4\sqrt{L}}{\eta rV}\right)\log^2(i+1)\right)\\
\leq&(\frac{n_k}{n})^{\delta}\left(\theta^*+\frac{\varepsilon_0}{V}\right)
+\frac{1}{n^{\delta}}S_{n_k}\left(\frac{2}{\eta}+\frac{4\sqrt{L}}{\eta rV}\right)\log^2n,
\end{align*}
where the first inequality follows from the fact that $\tilde\theta[n_k]<\theta^*+\varepsilon_0/V$.
Now, we take the $\limsup_{n\rightarrow\infty}$ from both sides and analyze each single term on the right hand side:
\begin{align*}
&1\geq\limsup_{n\rightarrow\infty}(\frac{n_k}{n})^{\delta}
 \geq\limsup_{k\rightarrow\infty}(\frac{n_k}{n_k+S_{n_k}})^{\delta}
 =\limsup_{k\rightarrow\infty}(\frac{1}{1+\frac{S_{n_k}}{n_k}})^{\delta}=1,~w.p.1,\\
&\limsup_{n\rightarrow\infty}\frac{S_{n_k}}{n^{\delta}}\left(\frac{2}{\eta}+\frac{4\sqrt{L}}{\eta rV}\right)\log^2n
\leq\limsup_{n\rightarrow\infty}\frac{S_{n_k}}{n_k^{1/3}}\cdot
\limsup_{n\rightarrow\infty}\frac{\left(\frac{2}{\eta}+\frac{4\sqrt{L}}{\eta rV}\right)\log^2n}{n^{\delta-1/3}}=0,~w.p.1,
\end{align*}
where we apply the convergence result \eqref{bound-on-returning} in the second line. Thus,
\[\limsup_{n\rightarrow\infty}\tilde{\theta}[n]\leq\theta^*+\frac{\varepsilon_0}{V},~w.p.1.\]
By Lemma \ref{truncation-lemma} we have $\limsup_{n\rightarrow\infty}\hat{\theta}[n]\leq\theta^*+\varepsilon_0/V$. Finally, by Lemma \ref{properties},  and the fact that $\theta^*+\varepsilon_0/V\in(0,\theta_{\max})$, we have  $\limsup_{n\rightarrow\infty}\theta[n]\leq\theta^*+\varepsilon_0/V$. Since this holds for any $\varepsilon_0>0$ small enough, let $\varepsilon_0\rightarrow0$ finishes the proof.
\Halmos
\endproof

\subsection{Finishing the proof of near optimality}
With the help of previous analysis on $\theta[n]$, we are ready to prove our main theorem, with the following lemma on strong law of large numbers for martingale difference sequences:
\begin{lemma}[Corollary 4.2 of \citet{neely2012stability}]\label{SLLN}
Let $\{\mathcal{F}_i\}_{i=0}^{\infty}$ be a filtration and let $\{X(i)\}_{i=0}^{\infty}$ be a real-valued random process such that $X(i)\in\mathcal{F}_{i+1},~\forall i$. Suppose there is a finite constant $C$ such that $\expect{X(i)|\mathcal{F}_i}\leq C,~\forall i$, and 
\[\sum_{i=1}^{\infty}\frac{\expect{X(i)^2}}{i^2}<\infty.\]
Then,
\[\limsup_{n\rightarrow\infty}\frac{1}{n}\sum_{i=0}^{n-1}X(i)\leq C,~~w.p.1.\]
\end{lemma}

\proof{Proof of Theorem \ref{theorem_average_converge}.~~}
Recall for any $n$, the empirical accumulation without ceil and floor function is
\[\hat{\theta}[n]=\frac{1}{n^{\delta}}\sum_{i=0}^{n-1}\left(y[i]-\theta[i]T[i]+\frac{1}{V}\sum_{l=1}^LQ_l[i](z_l[i]-c_lT[i])\right).\]
Dividing both sides by $\sum_{i=0}^{n-1}T[i]/n^{\delta}$ yields
\begin{align*}
\frac{\hat{\theta}[n]}{\frac{1}{n^{\delta}}\sum_{i=0}^{n-1}T[i]}
=&\frac{\sum_{i=0}^{n-1}\left(y[i]-\theta[i]T[i]+\frac{1}{V}\sum_{l=1}^LQ_l[i](z_l[i]-c_lT[i])\right)}{\sum_{i=0}^{n-1}T[i]}\\
=&\frac{\sum_{i=0}^{n-1}\left(y[i]+\frac{1}{V}\sum_{l=1}^LQ_l[i](z_l[i]-c_lT[i])\right)}{\sum_{i=0}^{n-1}T[i]}
    -\frac{\sum_{i=0}^{n-1}\theta[i]T[i]}{\sum_{i=0}^{n-1}T[i]}.
\end{align*}
Moving the last term to the left hand side and taking the $\limsup_{n\rightarrow\infty}$ from both sides gives
\begin{align*}
\limsup_{n\rightarrow\infty}\left(\frac{\hat{\theta}[n]}{\frac{1}{n^{\delta}}\sum_{i=0}^{n-1}T[i]}
+\frac{\sum_{i=0}^{n-1}\theta[i]T[i]}{\sum_{i=0}^{n-1}T[i]}\right)
\geq&\limsup_{n\rightarrow\infty}\frac{\sum_{i=0}^{n-1}y[i]}{\sum_{i=0}^{n-1}T[i]}
    +\frac{\sum_{i=0}^{n-1}\frac{1}{V}\sum_{l=1}^LQ_l[i](z_l[i]-c_lT[i])}{\sum_{i=0}^{n-1}T[i]}\\
\geq&\limsup_{n\rightarrow\infty}\frac{\sum_{i=0}^{n-1}y[i]}{\sum_{i=0}^{n-1}T[i]}
+\frac{1}{2}\frac{\|\mathbf{Q}[n]\|^2-\sum_{i=0}^{n-1}\sum_{l=1}^L(z_l[i]-c_lT[i])^2}{V\sum_{i=0}^{n-1}T[i]}\\
\geq&\limsup_{n\rightarrow\infty}\frac{\sum_{i=0}^{n-1}y[i]}{\sum_{i=0}^{n-1}T[i]}
-\frac{1}{2V}\limsup_{n\rightarrow\infty}\frac{1}{n}\sum_{i=0}^{n-1}K[i]^2,
\end{align*}
where the second inequality follows from inequality \eqref{dpp-relation} and telescoping sums, and the last inequality follows from $T[n]\geq1$, $\|\mathbf{Q}[n]\|^2\geq0$ and $K[i]=\sqrt{\sum_{l=1}^L(z_l[i]-c_lT[i])^2}$. Now we use Lemma \ref{SLLN} with $X(i)=K[i]^2$ to bound the second term. Since $K[i]$ is of exponential type by Assumption \ref{bounded-assumption}, we know that $\expect{K[i]^2|\mathcal{H}_n}\leq 2B^2/\eta^2$. Furthermore, $\expect{K[i]^4}\leq24B^4/\eta^4$. Thus,
\[\sum_{i=1}^{\infty}\frac{\expect{K[i]^4}}{i^2}<\infty.\]
Thus, all assumptions in Lemma \ref{SLLN} are satisfied and we conclude that
\[\limsup_{n\rightarrow\infty}\frac{1}{n}\sum_{i=0}^{n-1}K[i]^2\leq \frac{2B^2}{\eta^2},~w.p.1.\]
This implies,
\[\limsup_{n\rightarrow\infty}\left(\frac{\hat{\theta}[n]}{\frac{1}{n^{\delta}}\sum_{i=0}^{n-1}T[i]}
+\frac{\sum_{i=0}^{n-1}\theta[i]T[i]}{\sum_{i=0}^{n-1}T[i]}\right)
\geq\limsup_{n\rightarrow\infty}\frac{\sum_{i=0}^{n-1}y[i]}{\sum_{i=0}^{n-1}T[i]}-\frac{B^2}{\eta^2V}.\]
By Theorem \ref{theorem-asymptotic-upperbound}, $\hat{\theta}[n]$ is asymptotically upper bounded. Since $\delta<1$ and $T[n]\geq1$, it follows $\frac{1}{n^{\delta}}\sum_{i=0}^{n-1}T[i]=\mathcal{O}(n^{1-\delta})$, which goes to infinity as $n\rightarrow\infty$. Thus,
\[\limsup_{n\rightarrow\infty}\frac{\hat{\theta}[n]}{\frac{1}{n^{\delta}}\sum_{i=0}^{n-1}T[i]}\leq0,\]
and thus,
\[\limsup_{n\rightarrow\infty}\frac{\sum_{i=0}^{n-1}\theta[i]T[i]}{\sum_{i=0}^{n-1}T[i]}
\geq\limsup_{n\rightarrow\infty}\frac{\sum_{i=0}^{n-1}y[i]}{\sum_{i=0}^{n-1}T[i]}
-\frac{B^2}{\eta^2V}.\]
By Theorem \ref{theorem-asymptotic-upperbound} again, $\theta[n]$ is asymptotically upper bounded by $\theta^*$. Based on this result, it is easy to show the following
\[\limsup_{n\rightarrow\infty}\frac{\sum_{i=0}^{n-1}\theta[i]T[i]}{\sum_{i=0}^{n-1}T[i]}\leq\theta^*.\]
Thus, we finally get
\[\limsup_{n\rightarrow\infty}\frac{\sum_{i=0}^{n-1}y[i]}{\sum_{i=0}^{n-1}T[i]}\leq\theta^*+\frac{B^2}{\eta^2V},\]
finishing the proof.
\Halmos
\endproof

\section{Simulation experiments}\label{simulation}
In this section, we demonstrate the performance of our proposed algorithm through an application scenario on single user file downloading. We show that this problem can be formulated as a two state constrained online MDP and solved using our proposed algorithm.

Consider a slotted time system where $t\in\{0,1,2,\cdots\}$, and one user is repeatedly downloading files.
We use $F(t)\in\{0,1\}$ to denote the system file state at time slot $t$.  State ``1'' indicates there is an active file in the system for downloading and state ``0'' means there is no file and the system is idle.
Suppose the user can only download 1 file at each time, and the user cannot observe the file length. Each file contains an integer number of packets which is independent and geometrically distributed with expected length equal to 1.

During each time slot where there is an active file for downloading (i.e. $F(t)=1$), the user first observes the channel state
$\omega(t)$, which is the i.i.d. random variable taking values in $\Omega=\{0.2, 0.5, 0.8\}$ with equal probabilities, and delay penalty $s(t)$,  which is also an i.i.d. random variable taking values in $\{1,3,5\}$ with equal probability. Then, the user
makes a service action $\alpha(t)\in\mathcal{A}=\{0, 0.3, 0.6, 0.9\}$. The pair $(\omega(t),\alpha(t))$ affects the following quantities:
\begin{itemize}
\item The success probability of downloading a file at time $t$: $\phi(\alpha(t),\omega(t))\triangleq\alpha(t)\cdot\omega(t)$.
\item The resource consumption $p(\alpha(t))$ at time $t$. We assume $p(0)=0$, $p(0.3)=1$, $p(0.6)=2$ and $p(0.9)=4$.
\end{itemize}
After a file is downloaded, the system goes idle (i.e. $F(t)=0$) and stays there for a random amount of time that is independent and geometrically
distributed with mean equal to 2. The goal is to minimize the time average delay penalty subject to a resource constraint that the time average resource consumption cannot exceed 1.

In \citet{wei2015power}, a similar optimization problem is considered but without random events $\omega(t)$ and $s(t)$, which can be formulated as a two state constrained MDP. Here, using the same logic, we can formulate our optimization problem as a two state constrained online MDP.
Given $F(t)=1$, the file will finish its download at the end of this time slot with probability $\phi(\alpha(t),\omega(t))$. Thus, the transition probabilities out of state 1 are:
\begin{align*}
&Pr[F(t+1)=0|F(t)=1]=\phi(\alpha(t),\omega(t))\\
&Pr[F(t+1)=1|F(t)=1]=1-\phi(\alpha(t),\omega(t)),
\end{align*}
On the other hand, given $F(t)=0$, the system is idle and will transition to the active state in the next slot with probability $\lambda$:
\begin{align*}
&Pr[F(t+1)=1|F(t)=0]=\lambda\\
&Pr[F(t+1)=0|F(t)=0]=1-\lambda,
\end{align*}

Now, we characterize this online MDP through renewal frames and show that it can be solved using the proposed algorithm in Section \ref{formulation}. First, notice that the state ``1'' is recurrent under any action $\alpha(t)$. We denote $t_n$ as the $n$-th time slot when the system returns to state ``1''. Define the renewal frame as the time period between $t_n$ and $t_{n+1}$ with frame size
\[T[n]=t_{n+1}-t_n.\]
Furthermore, since the system does not have any control options in state ``0'', the controller makes exactly one decision during each frame and this decision is made at the beginning of each frame. Thus, we can write out the optimization problem as follows:
\begin{align*}
\min~~&\limsup_{N\rightarrow\infty}\frac{\sum_{n=0}^{N-1}\alpha(t_n)s(t_n)}{\sum_{n=0}^{N-1}T[n]}\\
s.t.~~&\limsup_{N\rightarrow\infty}\frac{\sum_{n=0}^{N-1}p(\alpha(t_n))}{\sum_{n=0}^{N-1}T[n]}\leq1,
~\alpha(t_n)\in\mathcal{A}.
\end{align*}
Subsequently, in order to apply our algorithm, we can define the virtual queue $Q[n]$ as $Q[0]=0$ with updating rule
\[Q[n+1]=\max\{Q[n]+p(\alpha(t_n))-T[n],0\}.\]
Notice that for any particular action $\alpha(t_n)\in\mathcal{A}$ and random event $\omega(t_n)\in\Omega$, we can always compute $\expect{T[n]}$ as
\begin{align*}
\expect{T[n]}&=1-\phi(\alpha(t_n),\omega(t_n))+\phi(\alpha(t_n),\omega(t_n))\left(1+\frac{1}{\lambda}\right)\\
&=1+2\alpha(t_n)\omega(t_n),
\end{align*}
where the second equality follows by substituting $\lambda=0.5$ and $\phi(\alpha(t_n),\omega(t_n))=\alpha(t_n)\omega(t_n)$. Thus, for each $\alpha(t_n)\in\mathcal{A}$, the expression \eqref{DPP} can be computed.

In each of the simulations, each data point is the time average of 2 million slots. We compare the performance of the proposed algorithm with the optimal randomized policy. The optimal policy is computed by formulating the MDP into a linear program with the knowledge of the distribution on $\omega(t)$ and $s(t)$. See \citet{Fo66} for details of this linear program formulation.

\begin{figure}[htbp]
   \centering
   \includegraphics[height=3.5in]{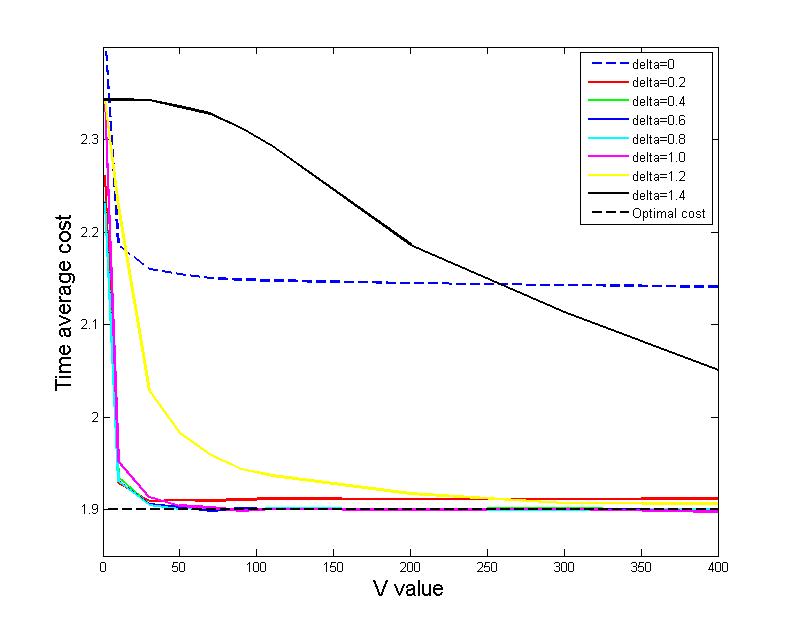} 
   \caption{Time average penalty versus tradeoff parameter V}
   \label{fig:Stupendous1}
\end{figure}

\begin{figure}[htbp]
   \centering
   \includegraphics[height=3.5in]{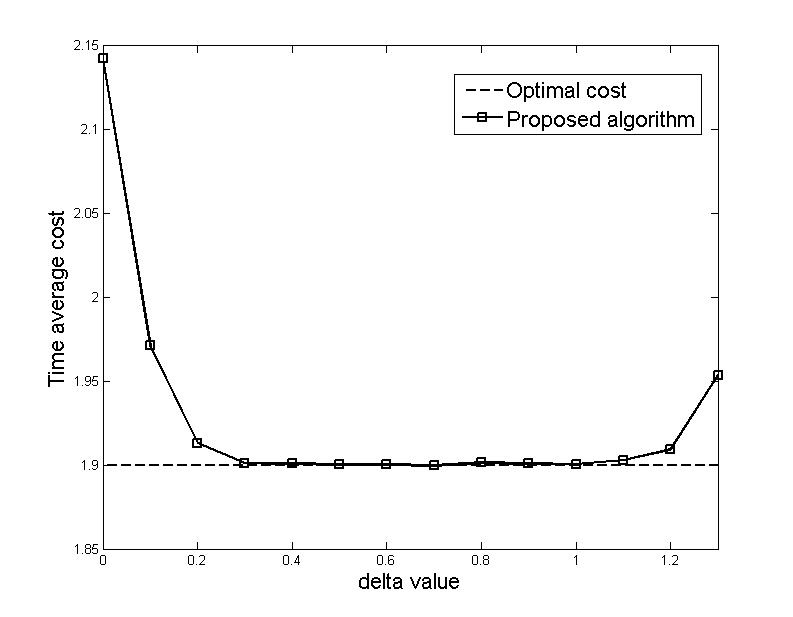} 
   \caption{Time average penalty versus $\delta$ parameter with fixed $V=300$.}
   \label{fig:Stupendous2}
\end{figure}

In Fig. \ref{fig:Stupendous1}, we plot the performance of our algorithm verses $V$ parameter for different $\delta$ value.
We see from the plots that as $V$ gets larger, the time averages approaches the optimal value and achieves a near optimal performance for $\delta$ roughly between $0.4$ and $1$. A more obvious relation between performance and $\delta$ value is shown in Fig. \ref{fig:Stupendous2}, where we fix $V=300$ and plot the performance of the algorithm verses $\delta$ value. It is clear from the plots that the algorithm fails whenever $\delta$ is too small ($\delta< 0.3$) or too big ($\delta>1$). This meets the statement of Theorem \ref{theorem_average_converge} that the algorithm works for $\delta\in(1/3,1)$.

\begin{figure}[htbp]
   \centering
   \includegraphics[height=3.5in]{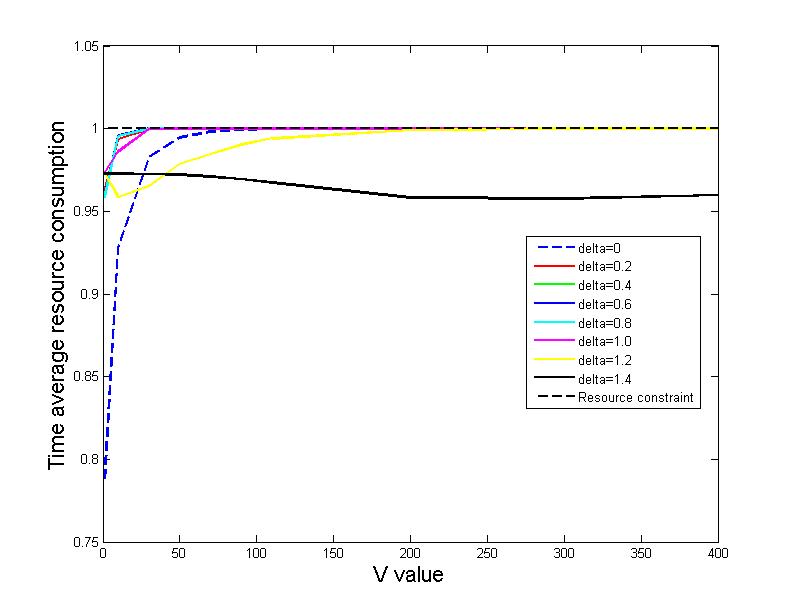} 
   \caption{Time average resource consumption versus tradeoff parameter $V$.}
   \label{fig:Stupendous3}
\end{figure}

\begin{figure}[htbp]
   \centering
   \includegraphics[height=3.5in]{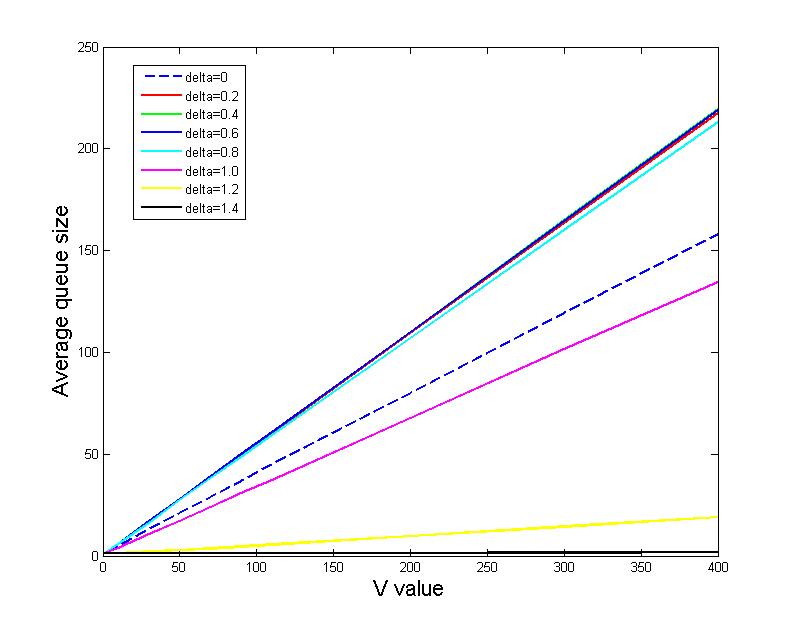} 
   \caption{Time average virtual queue size versus tradeoff parameter $V$.}
   \label{fig:Stupendous4}
\end{figure}

In Fig. \ref{fig:Stupendous3}, we plot the time average resource consumption verses $V$ value. We see from the plots that the algorithm is always feasible for different $V$'s and $\delta$'s, which meets the statement of Theorem \ref{feasibility}. Also, as $V$ gets larger, the constraint gap tends to be smaller.
In Fig. \ref{fig:Stupendous4}, we plot the average virtual queue size verses $V$ value. It shows that the average queue size gets larger as $V$ get larger. To see the implications,
recall from the proof of Theorem \ref{feasibility}, the inequality \eqref{inter-constraint-violation} implies that the virtual queue size $Q_l[N]$ affects the rate that the algorithm converges down to the feasible region. Thus, if the average virtual queue size is large, then, it takes longer for the algorithm to converge. This demonstrates that $V$ is indeed a trade-off parameter which trades the sub-optimality gap for the convergence rate.

\section{Conclusions}
This paper considers the constrained optimization over a renewal system with observed random events at the beginning of each renewal frame. We propose an online algorithm which does not need the knowledge of the distributions of random events. We prove that this proposed algorithm is feasible and achieves $O(\varepsilon)$ near optimality by constructing an exponential supermartingale. Simulation experiments demonstrates the near optimal performance of the proposed algorithm.

\bibliographystyle{informs2014}
\bibliography{constrained-online}

\begin{APPENDICES}

\section{Additional proofs}\label{sec:proof}

\proof{Proof of Lemma \ref{geometric-bound}.~~}
We begin by bounding the difference $\left|\|\mathbf{Q}[n+1]\|-\|\mathbf{Q}[n]\|\right|$ for any $n$:
\begin{align*}
\big|\|\mathbf{Q}[n+1]\|-\|\mathbf{Q}[n]\|\big|
\leq&\|\mathbf{Q}[n+1]- \mathbf{Q}[n]\|\\
=&\sqrt{\sum_{l=1}^L\big( \max\{Q_l[n] + z_l[n]-c_lT[n],~0\} - Q_l[n]\big)^2}\\
\leq&\sqrt{\sum_{l=1}^L(z_l[n]-c_lT[n])^2}=K[n],
\end{align*}
where the first inequality follows from triangle inequality and
the last inequality follows from the fact that for any $a,b\in\mathbb{R}$, $|\max\{a+b,0\}-a|\leq |b|$. 
Thus, it follows,
\begin{align*}
\left|\expect{\left.\|\mathbf{Q}[n+1]\|-\|\mathbf{Q}[n]\|\right|\mathcal{H}_n}\right|
\leq \expect{\left.K[n]\right|\mathcal{H}_n}\leq \frac{B}{\eta},
\end{align*}
which follows from Proposition \ref{prop-1}. Also, we have
\begin{align*}
\expect{\left.e^{r(\|\mathbf{Q}[n+1]\|-\|\mathbf{Q}[n]\|)}\right|\mathcal{H}_n}
\leq&\expect{\left.\exp\left(rK[n]\right)\right|\mathcal{H}_n}\\
\leq&\expect{\left.\exp\left(\eta K[n]\right)\right|\mathcal{H}_n}
\leq B\triangleq \Gamma
\end{align*}
where the second to last inequality follows by substituting the definition $r=\min\left\{\eta,\frac{\xi\eta^2}{4B}\right\}\leq\eta$ and the last inequality follows from Assumption \ref{bounded-assumption}.

Next, suppose $\|\mathbf{Q}[n]\|> \sigma\triangleq C_0V$. Then, since the proposed algorithm minimizes the term on the right hand side of \eqref{dpp-upperbound} over all possible decisions at frame $n$, it must achieve smaller value on that term compared to that of $\xi$-slackness policy $\alpha^{(\xi)}[n]$ specified in Assumption \ref{slack}. Formally, this is
\begin{align*}
&\expect{\left.\sum_{l=1}^LQ_l[n](z_l[n]-c_lT[n])+V(y[n]-\theta[n]T[n])~\right|~\mathcal{H}_n,\omega[n]}\\
\leq&\expect{\left.\sum_{l=1}^LQ_l[n](z_l^{(\xi)}[n]-c_lT^{(\xi)}[n])+V(y^{(\xi)}[n]-\theta[n]T^{(\xi)}[n])~\right|~\mathcal{H}_n,\omega[n]}.
\end{align*}
where we used the fact that $\theta[n]$ and $\mathbf{Q}[n]$ are in $\mathcal{H}_n$.
Substitute this bound into the right hand side of \eqref{dpp-upperbound}  and take expectation from both sides regarding $\omega[n]$ gives
\begin{align*}
&\expect{\Delta[n]+V(y[n]-\theta[n]T[n])~|~\mathcal{H}_n}\\
\leq&\expect{\left.\sum_{l=1}^LQ_l[n](z_l^{(\xi)}[n]-c_lT^{(\xi)}[n])+V(y^{(\xi)}[n]-\theta[n]T^{(\xi)}[n])~\right|~\mathcal{H}_n}
+B^2/\eta^2.
\end{align*}
Since $\Delta[n]=\frac12(\|\mathbf{Q}[n+1]\|^2-\|\mathbf{Q}[n]\|^2)$,
This implies
\begin{align*}
&\expect{\|\mathbf{Q}[n+1]\|^2-\|\mathbf{Q}[n]\|^2~|~\mathcal{H}_n}\\
\leq&2B^2/\eta^2+2\expect{\left.\sum_{l=1}^LQ_l[n](z_l^{(\xi)}[n]-c_lT^{(\xi)}[n])+V(y^{(\xi)}[n]-\theta[n]T^{(\xi)}[n])
-V(y[n]-\theta[n]T[n])
\right|\mathcal{H}_n}\\
\leq&2B^2/\eta^2+2\sum_{l=1}^LQ_l[n]\expect{\left.z_l^{(\xi)}[n]-c_lT^{(\xi)}[n]
\right|\mathcal{H}_n}+2V\frac{B+\theta_{\max}B}{\eta}\\
\leq&2B^2/\eta^2+2V\frac{B+\theta_{\max}B}{\eta}-2\xi\sum_{l=1}^LQ_l[n]\\
\leq&2B^2/\eta^2+2V\frac{B+\theta_{\max}B}{\eta}-2\xi\|\mathbf{Q}[n]\|,
\end{align*}
where the second inequality follows from applying Proposition \ref{prop-1} to bound $\expect{T[n]|\mathcal{H}_n}$ as well as the fact that $0<\theta[n]< \theta_{\max}$, and the third inequality
follows from the $\xi$-slackness property as well as the assumption that $z_l^{(\xi)}[n]$ is i.i.d. over slots and hence independent of $Q_l[n]$. This further implies
\begin{align*}
&\expect{\|\mathbf{Q}[n+1]\|^2~|~\mathcal{H}_n}\\
\leq&\|\mathbf{Q}[n]\|^2-2\xi\|\mathbf{Q}[n]\|+2B^2/\eta^2+2V\frac{B+\theta_{\max}B}{\eta}\\
=&\|\mathbf{Q}[n]\|^2-2\xi\|\mathbf{Q}[n]\|+2B^2/\eta^2+2V\frac{B+\theta_{\max}B}{\eta}-\frac{\xi^2}{4}+\frac{\xi^2}{4}\\
=&\|\mathbf{Q}[n]\|^2-2\xi\|\mathbf{Q}[n]\|+\frac{2B^2/\eta^2+2V\frac{B+\theta_{\max}B}{\eta}-\frac{\xi^2}{4}}{\xi}\cdot\xi+\frac{\xi^2}{4}\\
=&\|\mathbf{Q}[n]\|^2-2\xi\|\mathbf{Q}[n]\|+C_0V\cdot\xi+\frac{\xi^2}{4}\\
\leq&\|\mathbf{Q}[n]\|^2-\xi\|\mathbf{Q}[n]\|+\frac{\xi^2}{4}=\left(\|\mathbf{Q}[n]\|-\frac\xi2\right)^2,
\end{align*}
where we use the fact that $C_0=\frac{2B^2}{V\xi\eta^2}+\frac{2}{\xi}\frac{B+\theta_{\max}B}{\eta}-\frac{\xi}{4V}$ and also the assumption that $\|\mathbf{Q}[n]\|\geq C_0V$.
Now take the square root from both sides gives
\[\sqrt{\expect{\|\mathbf{Q}[n+1]\|^2~|~\mathcal{H}_n}}\leq\|\mathbf{Q}[n]\|-\frac\xi2.\]
By concavity of $\sqrt{x}$ function, we have $\expect{\left.\|\mathbf{Q}[n+1]\|~\right|~\mathcal{H}_n}\leq\sqrt{\expect{\|\mathbf{Q}[n+1]\|^2~|~\mathcal{H}_n}}$, thus,
\begin{equation}\label{pre-conclusion}
\expect{\left.\|\mathbf{Q}[n+1]\|~\right|~\mathcal{H}_n}\leq\|\mathbf{Q}[n]\|-\frac\xi2.
\end{equation}
Finally, we claim that this gives that under the condition $\|\mathbf{Q}[n]\|> \sigma\triangleq C_0V$,
\begin{equation}\label{conclusion}
\expect{\left.e^{r(\|\mathbf{Q}[n+1]\|-\|\mathbf{Q}[n]\|)}\right|\mathcal{H}_n}\leq\rho\triangleq
1-\frac{r\xi}{2}+\frac{2B}{\eta^2}r^2<1.
\end{equation}
To see this, we expand $\expect{\left.e^{r(\|\mathbf{Q}[n+1]\|-\|\mathbf{Q}[n]\|)}\right|\mathcal{H}_n}$ using Taylor series as follows:
\begin{align*}
&\expect{\left.e^{r(\|\mathbf{Q}[n+1]\|-\|\mathbf{Q}[n]\|)}\right|\mathcal{H}_n}\\
=&1+r\expect{\left.\|\mathbf{Q}[n+1]\|-\|\mathbf{Q}[n]\|\right|\mathcal{H}_n}
+r^2\sum_{k=2}^{\infty}\frac{r^{k-2}\expect{\left.(\|\mathbf{Q}[n+1]\|-\|\mathbf{Q}[n]\|)^k\right|\mathcal{H}_n}}{k!}\\
\leq&1- \frac{r\xi}{2}
+r^2\sum_{k=2}^{\infty}\frac{r^{k-2}\expect{\left.(\|\mathbf{Q}[n+1]\|-\|\mathbf{Q}[n]\|)^k\right|\mathcal{H}_n}}{k!}\\
\leq&1- \frac{r\xi}{2}
+r^2\sum_{k=2}^{\infty}\frac{\eta^{k-2}\expect{\left.(\|\mathbf{Q}[n+1]\|-\|\mathbf{Q}[n]\|)^k\right|\mathcal{H}_n}}{k!}\\
=&1- \frac{r\xi}{2}+r^2\frac{\left(\expect{\left.e^{\eta(\|\mathbf{Q}[n+1]\|-\|\mathbf{Q}[n]\|)}\right|\mathcal{H}_n}
-\eta\expect{\left.\|\mathbf{Q}[n+1]\|-\|\mathbf{Q}[n]\|\right|\mathcal{H}_n}-1\right)}{\eta^2}\\
\leq&1- \frac{r\xi}{2}+\frac{B+\eta\cdot\frac{B}{\eta}}{\eta^2}r^2\\
\leq&1-\frac{r\xi}{2}+\frac{2B}{\eta^2}r^2=\rho,
\end{align*}
where the first inequality follows from \eqref{pre-conclusion}, the second inequality follows from $r\leq\eta$, and the second
to last inequality follows from Proposition \ref{prop-1}.

Finally, notice that the above quadratic function on $r$ attains the minimum at the point
$r=\frac{\xi\eta^2}{4B}$ with value $1-\frac{\xi^2\eta^2}{8B}<1$, and this function is strictly decreasing when
$$r\in\left(0, \frac{\xi\eta^2}{4B}\right).$$
Thus, our choice of
$$r=\min\left\{\eta,\frac{\xi\eta^2}{4B}\right\}
\leq\frac{\xi\eta^2}{4B}$$
ensures that $\rho$ is strictly less than 1 and the proof is finished.
\Halmos\\
\endproof

\proof{Proof of Lemma \ref{properties}.~~}
If $\theta[n]=y$ for some $y\in[0,\theta_{\max}]$, then, $\hat\theta[n]$ falls into one of the following three cases:
\begin{itemize}
\item $\hat\theta[n]=y$.
\item $y=\theta_{\max}$ and $\hat\theta[n]>\theta_{\max}$.
\item $y=0$ and $\hat\theta[n]<0$.
\end{itemize}
Then, we prove the above four properties based on these three cases.

1) If $\theta[n]=y\geq x$ for some $y$, then, the first two cases immediately imply $\hat\theta[n]\geq x$. If $y=0$, then, we have $x\leq0$, which violates the assumption that $x\in(0,\theta_{\max})$. Thus, the third case is ruled out. On the other hand, if $\hat\theta[n]\geq x$, then, obviously, $\theta[n]\geq x$.

2) If $\theta[n]=y\leq x$ for some $y$, then the last two cases immediately imply $\hat\theta[n]\leq x$. If $y=\theta_{\max}$, then, we have $x\geq y_{\max}$, which violates the assumption that $x\in(0,\theta_{\max})$. Thus, the first case is ruled out. On the other hand, if $\hat\theta[n]\leq x$, then, obviously, $\theta[n]\leq x$.

3) If $\limsup_{n\rightarrow\infty}\theta[n]\leq x$, then, for any $\epsilon>0$ such that $x+\epsilon<y_{\max}$, there exists an $N$ large enough so that $\theta[n]\leq x+\epsilon,~\forall n\geq N$. Then, by property 2), $\hat\theta[n]\leq x+\epsilon,~\forall n\geq N$, which implies $\limsup_{n\rightarrow\infty}\hat\theta[n]\leq x+\epsilon$. Let $\epsilon\rightarrow0$ gives $\limsup_{n\rightarrow\infty}\hat\theta[n]\leq x$.  One the other hand, if $\limsup_{n\rightarrow\infty}\hat\theta[n]\leq x$, then, obviously, $\limsup_{n\rightarrow\infty}\theta[n]\leq x$.

4) If $\liminf_{n\rightarrow\infty}\theta[n]\geq x$, then, for any $\epsilon>0$ such that $x-\epsilon>0$ there exists an $N$ large enough so that $\theta[n]\geq x-\epsilon,~\forall n\geq N$. Then, by property 1), $\hat\theta[n]\geq x-\epsilon,~\forall n\geq N$, which implies $\limsup_{n\rightarrow\infty}\hat\theta[n]\leq x-\epsilon$. Let $\epsilon\rightarrow0$ gives $\limsup_{n\rightarrow\infty}\hat\theta[n]\geq x$. One the other hand, if $\limsup_{n\rightarrow\infty}\hat\theta[n]\leq x$, then, obviously, $\limsup_{n\rightarrow\infty}\theta[n]\leq x$.
\Halmos\\
\endproof

\proof{Proof of Lemma \ref{exp-supMG}.~~}
The proof is divided into two parts. The first part contains some technical preliminaries showing $G[n]$ is measurable respect to $\mathcal{H}_n,~\forall n\geq n_k+1$, and the second part contains computations to prove the supermartingale claim.
\begin{itemize}
\item \textit{Technical preliminaries:}
First of all, for any fixed $k$, since $n_k$ is a random variable on the integers, we need to justify that $\{\mathcal{H}_n\}_{n\geq n_k+1}$ is indeed a filtration.  
First, it is obvious that $n_k$ a valid stopping time, i.e. 
$$\{n_k\leq t\}\in\mathcal{H}_t,~\forall t\in\mathbb{N}.$$
Then, any $n=n_k+s$ with some constant $s\in\mathbb{N}^+$ is also a valid stopping time because
$$\{n\leq t\}=\{n_k\leq t-s\}\in\mathcal{H}_{(t-s)\vee0}\subseteq\mathcal{H}_t,~\forall t\in\mathbb{N},$$
where $a\vee b\triangleq\max\{a,b\}$. Thus, by definition of stopping time $\sigma$-algebra from \citet{Durrett}, we know that for any $n\geq n_k+1$, $\mathcal{H}_n$ can be written as the collection of all sets $A$ that have $A\cap\{n\leq t\}\in\mathcal{H}_t,~\forall t\in\mathbb{N}$\footnote{An intuitive interpretation is that when $n\leq t$, the set $A$ is contained in the information known until $t$.}. Now, pick $1\leq s_1\leq s_2$ as constants, and if a set $A\in\mathcal{H}_{n_k+s_1}$, then, 
$$A\cap\{n_k+s_2\leq t\}=A\cap\{n_k+s_1\leq t-(s_2-s_1)\}\in\mathcal{H}_{(t-(s_2-s_1))\vee 0}\subseteq\mathcal{H}_t.$$
Thus, $\mathcal{H}_{n_k+s_1}\subseteq\mathcal{H}_{n_k+s_2}$ and $\{\mathcal{H}_n\}_{n\geq n_k+1}$ is indeed a filtration.

Since $\tilde\theta[n_k+1]$ is determined by the realization up to frame $n_k$, it follows, for any $t\in\mathbb{N}^+$,
$$\{\tilde\theta[n_k+1]\geq\theta^*+\varepsilon_0/V\}\cap\{n_k+1\leq t\}
=\cup_{s=1}^t\{\tilde\theta[s]\geq\theta^*+\varepsilon_0/V\}\in\mathcal{H}_t,$$
which implies that $\{\tilde\theta[n_k+1]\geq\theta^*+\varepsilon_0/V\}\in\mathcal{H}_{n_k+1}$. Since 
$\{\mathcal{H}_n\}_{n\geq n_k+1}$ is a filtration, it follows $\{\tilde\theta[n_k+1]\geq\theta^*+\varepsilon_0/V\}\in\mathcal{H}_n$ for any $n\geq n_k+1$. By the same methodology, we can show that
$\{\tilde\theta[n]<\theta^*+\varepsilon_0/V\}\in\mathcal{H}_n,~\forall n\geq n_k+1$, which in turn implies, $\{S_{n_k}+n_k\leq n\}\in\mathcal{H}_n$ and $\{S_{n_k}\geq n-n_k+1\}\in\mathcal{H}_n$. 
Overall, the function $G[n]$ is measurable respect to $\mathcal{H}_n,~\forall n\geq n_k+1$.

\item \textit{Proof of supermartingale claim:} It is obvious that $|G[n]|<\infty$, thus, in order to prove $G[n]$ is a supermartingale, it is enough to show that
\begin{equation}\label{sup_MG_condition}
\expect{\left.G[n+1]-G[n]\right|\mathcal{H}_n}\leq0,~\forall n\geq n_k+1.
\end{equation}
First, on the set $\{S_{n_k}\leq n-n_k\}$, we have
\[\expect{\left.(G[n+1]-G[n])\mathbf{1}_{\{S_{n_k}+n_k\leq n\}}\right|\mathcal{H}_n}=\expect{\left.(G[n]-G[n])\mathbf{1}_{\{S_{n_k}+n_k\leq n\}}\right|\mathcal{H}_n}=0.\]
It is then sufficient to show the inequality \eqref{sup_MG_condition} holds on the set $\{S_{n_k}\geq n-n_k+1\}$. Since
\begin{align*}
&\expect{G[n+1]\mathbf{1}_{\{S_{n_k}\geq n-n_k+1\}}|\mathcal{H}_n}\\
=&\expect{\left.\frac{e^{\lambda_{n+1} F[(n+1)\wedge(n_k+S_{n_k})]}}{\prod_{i=n_k+1}^{(n+1)\wedge (n_k+S_{n_k})}\rho_i}~\right|~\mathcal{H}_{n}}\mathbf{1}_{\{\tilde\theta[n_k+1]\geq\theta^*+\varepsilon_0/V\}}\mathbf{1}_{\{S_{n_k}\geq n-n_k+1\}}\\
=&\expect{\left.\frac{e^{\lambda_{n+1} F[n+1]}}{\prod_{i=n_k+1}^{n+1}\rho_i}~\right|~\mathcal{H}_{n}}
 \mathbf{1}_{\{\tilde\theta[n_k+1]\geq\theta^*+\varepsilon_0/V\}}\mathbf{1}_{\{S_{n_k}\geq n-n_k+1\}}\\
=&\frac{e^{\lambda_{n+1} F[n]}}{\prod_{i=n_k+1}^{n}\rho_i}
  \expect{\left.\frac{e^{\lambda_{n+1} (F[n+1]-F[n])}}{\rho_{n+1}}~\right|~\mathcal{H}_{n}}
  \mathbf{1}_{\{\tilde\theta[n_k+1]\geq\theta^*+\varepsilon_0/V\}}\mathbf{1}_{\{S_{n_k}\geq n-n_k+1\}}\\
\leq&\frac{e^{\lambda_{n} F[n]}}{\prod_{i=n_k+1}^{n}\rho_i}
  \expect{\left.\frac{e^{\lambda_{n+1} (F[n+1]-F[n])}}{\rho_{n+1}}~\right|~\mathcal{H}_{n}}
  \mathbf{1}_{\{\tilde\theta[n_k+1]\geq\theta^*+\varepsilon_0/V\}}\mathbf{1}_{\{S_{n_k}\geq n-n_k+1\}}\\
=&G[n]\expect{\left.\frac{e^{\lambda_{n+1} (F[n+1]-F[n])}}{\rho_{n+1}}~\right|~\mathcal{H}_{n}}
\mathbf{1}_{\{\tilde\theta[n_k+1]\geq\theta^*+\varepsilon_0/V\}}\mathbf{1}_{\{S_{n_k}\geq n-n_k+1\}},
\end{align*}
where $\mathbf{1}_{\{\tilde\theta[n_k+1]\geq\theta^*+\varepsilon_0/V\}}$ and $\mathbf{1}_{\{S_{n_k}\geq n-n_k+1\}}$ can be moved out of the expectation because $\{\tilde\theta[n_k+1]\geq\theta^*+\varepsilon_0/V\}\in\mathcal{H}_{n}$ and
$\{S_{n_k}\geq n-n_k+1\}\in\mathcal{H}_{n}$,
and the only inequality follows from the following argument: On the set $\{S_{n_k}\geq n-n_k+1\}$, $\{\tilde\theta[n]\geq\theta^*+\varepsilon_0/V\}$, thus, by Lemma \ref{comparison-lemma}, $F[n]\geq0$ and
using the fact $\lambda_n>\lambda_{n+1}$, we have $\lambda_{n+1} F[n]\leq\lambda_{n} F[n]$.
Thus, it is sufficient to show that on the set $\{S_{n_k}\geq n-n_k+1\}\cap\{\tilde\theta[n_k+1]\geq\theta^*+\varepsilon_0/V\}$, we have
\[\expect{\left.\frac{e^{\lambda_{n+1} (F[n+1]-F[n])}}{\rho_{n+1}}~\right|~\mathcal{H}_{n}}\leq1.\]
By Taylor expansion, we have
\begin{align*}
&\expect{\left.e^{\lambda_{n+1} (F[n+1]-F[n])}~\right|~\mathcal{H}_{n}}\\
=&1+\lambda_{n+1}\expect{F[n+1]-F[n]~|~\mathcal{H}_{n}}+\sum_{k=2}^{\infty}\frac{\lambda_{n+1}^k}{k!}\expect{(F[n+1]-F[n])^k~|~\mathcal{H}_{n}}\\
=&1+\lambda_{n+1}\expect{F[n+1]-F[n]~|~\mathcal{H}_{n}}
+\lambda_{n+1}^2\sum_{k=2}^{\infty}\frac{\lambda_{n+1}^{k-2}}{k!}\expect{(F[n+1]-F[n])^k~|~\mathcal{H}_{n}}\\
\leq&1-\frac{\lambda_{n+1}\varepsilon_0}{V}
+\lambda_{n+1}^2\sum_{k=2}^{\infty}\frac{\lambda_{n+1}^{k-2}}{k!}\expect{(F[n+1]-F[n])^k~|~\mathcal{H}_{n}},
\end{align*}
where the last inequality comes from the following argument: On the set $\{S_{n_k}\geq n-n_k+1\}$, $\tilde{\theta}[n_k+1]\geq\theta^*+\varepsilon_0/V$, thus,
by the definition of $\tilde\theta[n]$, we have
$\hat{\theta}[n]\geq\tilde\theta[n]\geq\theta^*+\varepsilon_0/V$, and Lemma \ref{properties} gives
$\theta[n]\geq\theta^*+\varepsilon_0/V$, then, by Lemma \ref{key-feature}, we have
\[\expect{F[n+1]-F[n]~|~\mathcal{H}_{n}}\leq-\frac{\varepsilon_0}{V}.\]

Now, by the assumption that $V\geq\frac{\varepsilon_0\eta}{4\log^22}-\frac{2\sqrt{L}}{r}$, we have
$\lambda_{n+1}\leq\frac{1}{\left(\frac{2}{\eta}+\frac{4\sqrt{L}}{\eta rV}\right)\log^2(n+1)}$, which follows from simple algebraic manipulations. Using the fact that $|F[n+1]-F[n]|\leq \left(\frac{2}{\eta}+\frac{4\sqrt{L}}{\eta rV}\right)\log^2(n+1)$, we have
\begin{align*}
&\expect{\left.e^{\lambda_{n+1} (F[n+1]-F[n])}~\right|~\mathcal{H}_{n}}\\
\leq&1-\frac{\lambda_{n+1}\epsilon_0}{V}+\lambda_{n+1}^2\sum_{k=2}^{\infty}
\frac{\left(\frac{1}{\left(\frac{2}{\eta}+\frac{4\sqrt{L}}{\eta rV}\right)\log^2(n+1)}\right)^{k-2}}{k!}
\expect{\left.\left(\left(\frac{2}{\eta}+\frac{4\sqrt{L}}{\eta rV}\right)\log^2(n+1)\right)^k~\right|~\mathcal{H}_{n}}\\
=&1-\frac{\lambda_{n+1}\epsilon_0}{V}+\lambda_{n+1}^2\sum_{k=2}^{\infty}\frac{1}{k!}
\left(\left(\frac{2}{\eta}+\frac{4\sqrt{L}}{\eta rV}\right)\log^2(n+1)\right)^2\\
\leq&1-\frac{\lambda_{n+1}\epsilon_0}{V}+\lambda_{n+1}^2e\left(\frac{2}{\eta}+\frac{4\sqrt{L}}{\eta rV}\right)^2\log^4(n+1)=\rho_{n+1},
\end{align*}
where the final inequality follows by completing the third term back to Taylor series which is equal to $e$.
Overall, the inequality \eqref{sup_MG_condition} holds and $G[n]$ is a supermartingale. \Halmos
\end{itemize}
.\\
\endproof

\section{Computation of Asymptotics}\label{computation}
In this appendix, we show that there exists a constant $C$ such that
\[\sum_{m=1}^\infty(m+1)^3\prod_{i=n_k+1}^{n_k+m}\rho_i\leq C(n_k+2)^{4\beta}.\]
We first bound $\rho_i$. Let $C_1=\frac{96V^2e\left(\frac{2}{\eta}+\frac{4\sqrt{L}}{\eta rV}\right)^2}{\varepsilon_0^2\beta^4}$, then,
\begin{align*}
\rho_i = &1-\frac{\varepsilon_0^2}{4V^2e\left(\frac{2}{\eta}+\frac{4\sqrt{L}}{\eta rV}\right)^2\log^4(i+1)} \\
&= 1 - \frac{1}{C_1\frac{\beta^4}{24}\log^4(i+1)}\\
&< 1 - \frac{1}{C_1(i + 1)^{\beta}},
\end{align*}
where we used the fact that $\frac{\beta^4}{24}\log^4(i+1)<(i+1)^\beta,~\forall \beta>0, i\geq0$. Next, to bound $\prod_{i=n_k+1}^{n_k+m}\rho_i$, we take the logarithm:
\begin{align*}
\log\left(\prod_{i=n_k+1}^{n_k+m}\rho_i\right)=&\sum_{i=n_k+1}^{n_k+m}\log\rho_i\\
=&\sum_{i=n_k+1}^{n_k+m}\log\left(1 - \frac{1}{C_1(i + 1)^{\beta}}\right)\\
\leq&-\sum_{i=n_k+1}^{n_k+m}\frac{1}{C_1(i + 1)^{\beta}}\\
\leq&-\frac{1}{C_1}\int_{n_k+2}^{n_k+m+1}\frac{1}{x^{\beta}}dx.
\end{align*}
where the first inequality follows from the first order Taylor expansion. Since $\beta<1$, we compute the integral, which gives
\[-\frac{1}{C_1}\int_{n_k+2}^{n_k+m+1}\frac{1}{x^{\beta}}dx
=-\frac{1}{C_1(1-2\beta)}\left((n_k+m+1)^{1-\beta}-(n_k+2)^{1-\beta}\right).\]
Thus,
\begin{align*}
&\sum_{m=1}^\infty(m+1)^3\prod_{i=n_k+1}^{n_k+m}\rho_i\\
\leq&\sum_{m=1}^\infty(m+1)^3e^{-\frac{1}{C_1(1-\beta)}\left((n_k+m+1)^{1-\beta}-(n_k+2)^{1-\beta}\right)}\\
\leq&\int_0^{\infty}(x+2)^3e^{-\frac{1}{C_1(1-\beta)}\left((x+n_k+2)^{1-\beta}-(n_k+2)^{1-\beta}\right)}dx
+(3C_1(1-\beta))^4,
\end{align*}
where the last inequality follows from the fact that the integrand is monotonically decreasing when $x>3C_1(1-\beta)$, thus, the integral dominates the sum on the tail $x>3C_1(1-\beta)$. For the part where $x\leq 3C_1(1-\beta)$, the maximum of the integrand is bounded by $(3C_1(1-\beta))^3$. Thus, the total difference of such approximation is bounded by $(3C_1(1-\beta))^4$.
Then, we try to estimate the integral. Notice that
\[\frac{d}{dx}e^{-\frac{1}{C_1(1-\beta)}(x+n_k+2)^{1-\beta}}
=-\frac{1}{C_1}e^{-\frac{1}{C_1(1-\beta)}(x+n_k+2)^{1-\beta}}(x+n_k+2)^{-\beta},\]
we do integration-by-parts, which gives
\begin{align*}
&\int_0^{\infty}(x+2)^3e^{-\frac{1}{C_1(1-\beta)}\left((x+n_k+2)^{1-\beta}-(n_k+2)^{1-\beta}\right)}dx\\
=&\int_0^{\infty}(x+2)^3(x+n_k+2)^{\beta}(x+n_k+2)^{-\beta}e^{-\frac{1}{C_1(1-\beta)}(x+n_k+2)^{1-\beta}}dx
\cdot e^{\frac{1}{C_1(1-\beta)}(n_k+2)^{1-\beta}}\\
=&8C_1(n_k+2)^{\beta}+\int_{0}^{\infty}C_1\left(3(x+2)^2(x+n_k+2)^{\beta}+\beta(x+2)^3(x+n_k+2)^{\beta-1}\right)
e^{-\frac{1}{C_1(1-\beta)}\left((x+n_k+2)^{1-\beta}-(n_k+2)^{1-\beta}\right)}dx.
\end{align*}
Since $5\beta\leq1$ and $n_k\geq1$, we have $x+n_k+2\geq x+2$, which implies $(x+2)^3(x+n_k+2)^{\beta-1}\leq(x+2)^2(x+n_k+2)^{\beta}$, thus,
\begin{align*}
&\int_0^{\infty}(x+2)^3e^{-\frac{1}{C_1(1-\beta)}\left((x+n_k+2)^{1-\beta}-(n_k+2)^{1-\beta}\right)}dx\\
\leq&8C_1(n_k+2)^{\beta}+\int_{0}^{\infty}4C_1(x+2)^2(x+n_k+2)^{\beta}e^{-\frac{1}{C_1(1-\beta)}\left((x+n_k+2)^{1-\beta}-(n_k+2)^{1-\beta}\right)}dx.
\end{align*}
Repeat above procedure 3 more times, we have
\begin{align*}
&\int_0^{\infty}(x+2)^3e^{-\frac{1}{C_1(1-\beta)}\left((x+n_k+2)^{1-\beta}-(n_k+2)^{1-\beta}\right)}dx\\
\leq&8C_1(n_k+2)^{\beta}+16C_1^2(n_k+2)^{2\beta}+24C_1^3(n_k+2)^{3\beta}+24C_1^4(n_k+2)^{4\beta}\\
&+\int_0^{\infty}24C_1^4(x+n_k+2)^{4\beta-1}e^{-\frac{1}{C_1(1-\beta)}\left((x+n_k+2)^{1-\beta}-(n_k+2)^{1-\beta}\right)}dx\\
\leq&8C_1(n_k+2)^{\beta}+16C_1^2(n_k+2)^{2\beta}+24C_1^3(n_k+2)^{3\beta}+24C_1^4(n_k+2)^{4\beta}
+24C_1^5\leq C(n_k+2)^{4\beta},
\end{align*}
for some $C$ on the order of $C_1^5$ (which is $\mathcal{O}\left(V^{10}\beta^{-20}\varepsilon_0^{-10}\right)$),
where the second to last inequality follows from $4\beta-1\leq-\beta$ and thus, we replace $(x+n_k+2)^{4\beta-1}$ with $(x+n_k+2)^{-\beta}$ and do a direct integration. Overall, we proved the claim.

\end{APPENDICES}

\end{document}